\documentclass[reqno,12pt]{amsart}  

 \usepackage{gensymb}
 
 \usepackage{amsmath}
\usepackage{amssymb} 
\usepackage{mathtools}   
\usepackage{graphicx}
\usepackage{graphics}  
\usepackage{color}
\usepackage{verbatim} 
\usepackage[normalem]{ulem} 
\usepackage[mathscr]{eucal}
\usepackage{upgreek}
\usepackage{enumerate} 
 \usepackage[english]{babel} 
 \usepackage{float}
 \usepackage{subfigure}

\usepackage{babel,blindtext} 
 
\usepackage[linecolor=red,backgroundcolor=red!25,bordercolor=red,textsize=tiny]{todonotes} 
\usepackage{bm}

\floatstyle{boxed}
\restylefloat{figure}

\usepackage[titletoc]{appendix}
\numberwithin{equation}{section}

\usepackage[refpage,noprefix]{nomencl}
\usepackage{nomencl} 

\makenomenclature

\newtheorem{theorem}{Theorem}[section] 
\newtheorem{lemma}[theorem]{Lemma} 
\newtheorem{proposition}[theorem] {Proposition} 
\newtheorem{cor}[theorem]  {Corollary} 
\newtheorem{remark}[theorem]  {Remark} 
\newtheorem{definition}[theorem] {Definition}

\theoremstyle{definition}

\DeclareMathAlphabet{\mathpzc}{OT1}{pzc}{m}{it}

\DeclarePairedDelimiter{\abs}{\lvert}{\rvert}

\newcommand{\bF} {\boldsymbol{F}}

\newcommand{\bO}{\boldsymbol{\Omega}}

\newcommand{\bx} {\boldsymbol{x}}

\newcommand{\bGamma} {\boldsymbol{\Gamma}}

\newcommand{\bo}{\boldsymbol{\omega}}
\newcommand{\btau}{\boldsymbol{\tau}}

\newcommand{\bzeta}{\boldsymbol{\zeta}}

\newcommand{\bDel}{\boldsymbol{\Del}}

\newcommand{\bPi}{{\boldsymbol{\Pi}}} 
\newcommand{\bHcal}{\boldsymbol{\Hcal}}
\newcommand{\bFcal}{\boldsymbol{\Fcal}}

\renewcommand{\L} {\Lambda} 
\renewcommand{\O} {\Omega}

\newcommand{\eps}{\varepsilon}

\newcommand{\om}{\omega}

\newfam\Bbbfam 
\font\tenBbb=msbm10 
\font\sevenBbb=msbm7 
\font\fiveBbb=msbm5 
\textfont\Bbbfam=\tenBbb 
\scriptfont\Bbbfam=\sevenBbb 
\scriptscriptfont\Bbbfam=\fiveBbb

\newcommand{\R}     {\mathbb{R}} 
\newcommand{\Z}     {\mathbb{Z}} 
\newcommand{\N}     {\mathbb{N}} 
\renewcommand{\P}   {\mathbb{P}} 
 
\newcommand{\E}     {\mathbb{E}}

\newcommand{\1}{{\mathchoice {1\mskip-4mu\mathrm l}      
{1\mskip-4mu\mathrm l} 
{1\mskip-4.5mu\mathrm l} {1\mskip-5mu\mathrm l}}} 
\newcommand{\ssup}[1] {{\scriptscriptstyle{({#1}})}} 
\newtheoremstyle{thm}{2ex}{2ex}{\itshape\rmfamily}{} 
{\bfseries\rmfamily}{}{1.7ex}{} 
 \newtheoremstyle{rem}{1.3ex}{1.3ex}{\rmfamily}{} 
{\itshape\rmfamily}{}{1.5ex}{} 
 \newenvironment{proofsect}[1] 
{\vskip0.1cm\noindent{\scshape #1.}\hskip0.5cm}

\newcommand{\Acal}   {{\mathcal A }}

\newcommand{\Ecal}   {{\mathcal E }} 
\newcommand{\Fcal}   {{\mathcal F }} 
 
\newcommand{\Hcal}   {{\mathcal H }}

\newcommand{\Mcal}   {{\mathcal M }}

\newcommand{\Lscr} {\mathscr{L}}

\newcommand{\Escr}{{\mathscr{E}}}

\newcommand{\Gscr}{\mathscr{G}}

\newcommand{\Asf}{{\mathsf{A}}}
\newcommand{\Csf}{{\mathsf{C}}}
\newcommand{\Zsf}{{\mathsf{Z}}}
\newcommand{\Psf}{{\mathsf{P}}}

 \newcommand{\ex}{{\rm e}} 
 
\renewcommand{\d}{{\rm d}}

\newcommand{\Leb}{{\rm Leb}}

\newcommand{\Del}{{\operatorname{\sf Del}}}
\newcommand{\Bb}{{\operatorname{\sf B}}}
\newcommand{\Vor}{{\operatorname {\sf Vor}}}

\newcommand{\Exp}{\mathscr{E}\kern-0.2mm{\operatorname{xp}}}
\newcommand{\Log}{\mathscr{L}\kern-0.2mm{\operatorname{og}}}

\newcommand{\pr}{{\operatorname {pr}}}
\renewcommand{\emptyset} {\varnothing}

\DeclareMathOperator*{\essinf}{ess\,inf} 
 
\makeindex

\begin{document}

\title[\hfill Phase transitions\hfill]
{Phase transition for Gibbs Delaunay Tessellations with geometric hardcore conditions}


\author{Stefan Adams and  Shannon Horrigan}
\address{Mathematics Institute, University of Warwick, Coventry CV4 7AL, United Kingdom}
\email{S.Adams@warwick.ac.uk}

\thanks{}
  

 


\begin{abstract}
In this paper, we prove the existence of infinite Gibbs Delaunay Potts tessellations for marked particle configurations. The particle systems has two types of interaction, a so-called \emph{background potential} ensures that small and large triangles are excluded in the Delaunay tessellation, and is similar to the so-called hardcore potential introduced in \cite{Der08}. Particles carry one of $q$ separate marks. Our main result is that for large activities and high \emph{type interaction} strength the model has at least $ q$ distinct translation-invariant Gibbs Delaunay Potts tessellations. The main technique is a coarse-graining procedure using the scales in the system followed by comparison with site percolation on $ \Z^2 $. 
\end{abstract}

\maketitle


\section{Introduction}

\subsection{Background} The goal of equilibrium statistical mechanics is to explain the macroscopic behaviour of physical systems in thermodynamic equilibrium in terms of the interactions between its microscopic elements. The concept of a Gibbs measure was introduced independently by Dobrushin \cite{Dob68} and Lanford and Ruelle \cite{LR69} as a mathematical description of an equilibrium state of a system consisting of a large number of interacting components. Many physical systems exhibit phase transitions in which the system moves from one equilibrium state to another, such as the magnetisation of a ferromagnetic metal such as iron, or the transition of a real gas from liquid to vapour. This phenomenon should manifest in our mathematical model through the existence of multiple Gibbs measures. If this is the case, we say that the model exhibits a phase transition.
The first model of a system of locally interacting particles in which phase transition was shown to occur is called the Ising model. It was introduced in 1920 by Wilhelm Lenz, with the hope of obtaining an understanding of ferromagnetic behaviour. In the Ising model each vertex of a lattice is envisaged as a particle and randomly assigned a `spin'  value of $ -1 $ or $1 $ representing its magnetic moment, with interactions between particles which favour the agreement of neighbouring spins. In dimensions $ d \ge 2 $ there is a critical temperature below which the interactions are strong enough that the magnetic moments can align and one of the two spin types dominates the other.
The $q$-state Potts model for $q\in\N $ is a generalisation of the Ising model where the spin value assigned to a particle (now called its `type') can take any value in the set $\{1,\ldots,q\} $.  This model exhibits a similar break of symmetry as the Ising model; there is a critical temperature below which particles of one type dominate the others. This results in the existence of q distinct Gibbs measures. This result for the Potts model was proven using the Fortuin-Kasteleyn representation in which the q-state Potts model is coupled with a bond percolation model on $\Z^d $ called the random cluster model with the same parameter $q$. The random cluster model is a generalisation of the Bernoulli bond percolation model with an additional weighting applied to each configuration depending on the parameter $q$  and the number of connected components the configuration has. Using this representation it is shown that percolation in the random cluster model coincides with the existence of multiple Gibbs measures in the Potts model.
The study of phase transitions is one of the main subjects of statistical mechanics, but even so, most models which are known to exhibit a phase transition are discrete like the two discussed above. We are interested here in the existence of phase transitions in continuum particle systems. Notable examples of phase transitions in the continuous setting include a gas-liquid phase transition \cite{LMP99} and the spontaneous breaking of rotational symmetry in a simple model of a two dimensional crystal without defects for small temperatures \cite{MR09}. Phase transitions are also established in \cite{GH96} for a class of continuum multi-type particle systems (continuum Potts models) in $\R^d$  for $d\ge 2$ with a finite range repulsive pair interaction between particles of different types. This class includes the Widom-Rowlinson model. The approach taken by the authors involves a random cluster representation analogous to the Fortuin-Kasteleyn representation in the discrete case. Percolation in the relevant random cluster model then implies the existence of multiple Gibbs measures. In \cite{AE16} the particle interaction from \cite{GH96} between unlike particles was replaced by an interaction based on the geometry of the Delaunay graph with a hard-core background potential depending on the lengths of the Delaunay edges. Both edge and triangle-dependent interactions were considered, and a random cluster representation was used in the proof. In a continuation \cite{AE19}, the authors also obtained a phase transition result with no background potential and finite range interaction depending on the lengths of the Delaunay edges.

\subsection{Outline}
In this manuscript we investigate a model in which there is a \textit{type potential} and a type independent \textit{background potential}, both of which depend on the geometry of each Delaunay triangle. The type potential depends on the triangle's area and has the effect that triangles with smaller areas are more likely to have all of their vertices belonging to the same type. The background potential is  equivalent to that in \cite{Der08}, and introduces hardcore constraints excluding both small and large triangles. In Section~\ref{sec-existence} we show in Theorem~\ref{Propexist} that infinite volume Gibbs Delaunay Potts tessellations exist using methods established in \cite{DDG12}. Our main result (Theorem~\ref{THM-main}) in Section~\ref{sec-phase} is that a phase transition occurs in this model for certain choices of the model parameters. We accomplish the proof using a random cluster representation as in \cite{AE16,AE19}. 

We introduce definitions and notations concerning Delaunay tessellations and our model in Section~\ref{sec-def} below. The random cluster representation we introduce in Section~\ref{sec-DelCluster} couples the finite volume Gibbs distribution in $\L\subset\R^2 $ (which we call the Delaunay continuum Potts distribution in $\L$) with a second measure known as the Delaunay random cluster distribution in $\L$. This coupling is restricted to the case when the former has a boundary condition made up of particles of the same type, which corresponds to a Delaunay random cluster distribution with a  `wired' boundary condition in which all hyperedges sufficiently far away from $\L$ are open.  

Finally, our \emph{background potential} leads to some regular structure of the random Delaunay tessellations and one may wonder if the \emph{type interaction}  alone can trigger a phase transition. This is a challenging and open question which we hope to address in the near future. The idea is that the \emph{type potential} increases the chances that a small triangle is open in the  Delaunay random cluster measure leading to some control on the number of connected components in the percolation model. On the other hand the range of the \emph{type potential} is chosen such that triangle with large areas are open once the inverse temperature is large enough.
\subsection{Remarks on Delaunay tessellations}
 Our work is partly motivated by results in \cite{Der08} where for the first time a so-called double hardcore potential has been introduced and the existence of infinite volume Gibbs Delaunay tessellations has been proved. The advantage of the hardcore potential is that our random Delaunay tessellations possess sufficient regular structure to make them accessible for applying a coarse-graining procedure and ultimately proving phase transitions. 
 
  There are differences between geometric models on the Delaunay hypergraph structure and point particle  models such as the Widom-Rowlinson model. The first is that edges and triangles in the Delaunay hypergraph are each proportional in number to the number of points  in the configuration. However, in the case of the complete hypergraph the number of edges is proportional to the number of points  squared and the number of triangles is proportional to the number of points  cubed. Secondly, in complete graphs of all classical models, the neighbourhood of a given point depends only on the distance between points and so the number of neighbours increases with the activity parameter $ z $ of the underlying point process. This means that the system will become strongly connected for high values of $z$. This is not the case for the Delaunay hypergraphs which exhibit a self-similar property. Essentially, as the activity parameter $z$ increases, the expected number of neighbours to a given point in the Delaunay hypergraph remains the same, see \cite{M94}. Therefore, in order to keep a strong connectivity, we use a type interaction between points of Delaunay edges with differing marks. Finally, and perhaps most importantly, is the question of additivity. Namely, suppose we have an existing point configuration $\om$ and we want to add a new particle $x$ to it. It is well-know that  classical many-body interactions are additive, see \cite{AE16} for details and references. On the other hand, in the Delaunay framework, the introduction of a new point to an existing configuration not only creates new edges and triangles, but destroys some too. The Delaunay interactions are therefore not additive, and for this reason, attractive and repulsive interactions are indistinct. In the case of a hard exclusion interaction, we arrive at the possibility that a configuration $ \om $ is excluded, but for some $x$, $ \om\cup x $ is not. This is called the non-hereditary property \cite{DG09}, which seems to rule out using techniques such as stochastic comparisons of point processes \cite{GK97}.

\section{Definitions and Notations}\label{sec-def}
\subsection{Setup}

We consider systems of particles  in $ \R^2 $, both in the case where the particles are described by their spatial location only, and where the particles possess a \emph{mark} describing their type or internal degree of freedom. The \emph{mark space} is the set $ Q=\{1,\ldots,q\},$ where $q\in\N, q\ge 2 $. Each marked point lies in the set $ \R^2\times Q $, and each marked configuration $ \bo $ is a countable subset of $ \R^2\times Q $ having a locally finite projection onto $ \R^2 $. We denote by 
$$ 
\bO= \{\bo\subset\R^2\times Q\colon \omega\in\Omega\}
$$ the set of all marked configurations with locally finite projection onto $ \R^2 $. Here, $\omega$ is the projection of $\bo$ onto $\R^2$ and $\O$ denotes the set of locally finite subsets of $ \R^2 $. We will sometimes identify $ \bo $ with a vector $ \bo=(\omega^{\ssup{1}},\ldots,\om^{\ssup{q}}) $ of pairwise disjoint locally finite sets $ \om^{\ssup{1}},\ldots,\om^{\ssup{q}} $ in $ \R^2 $. Each $ \bo\in\bO $ has an associated \emph{mark function} $ \sigma_\omega\in  Q^\omega $ which retrieves the mark of a point given its position, i.e., $ \sigma_\omega(x) = i $ if $  (x,i)\in\bo $. A marked configuration $\bo$ can therefore also be represented as a pair $(\omega, \sigma_{\omega}).$

For each measurable set $B$  in $ \R^2\times Q  $ the counting variable $ N(B)\colon\bo\to\bo(B) $ on $ \bO $ gives the number of marked particles such that  the pair (position, mark) belongs to $B$. We equip the space $ \bO $ with the $\sigma$-algebra $\boldsymbol{\Fcal} $ generated by the counting variables $N(B) $ and the space  $ \O $ of locally finite configurations with the $\sigma$-algebra $ \Fcal $ generated by the counting variables $ N_\Delta: \om \mapsto \#\{\om\cap\Delta\} $ for $ \Delta\Subset\R^2 $ where we write $ \Delta\Subset\R^2 $ for any bounded $ \Delta\subset\R^2 $. As usual, we take as the reference measure on $ (\bO,\boldsymbol{\Fcal}) $ the marked Poisson point process $ \boldsymbol{\Pi}^z $ with intensity measure $ z\Leb\otimes\mathsf{U} $  where $ z>0 $ is an arbitrary activity, $ \Leb $ is the Lebesgue measure in $ \R^2 $, and $ \mathsf{U} $ is the uniform probability measure on $ Q  $.

For each $ \L\subset\R^2 $ we write $ \bO_\L=\{\bo\in\bO\colon \bo\subset\L\times  Q\} $  for the set of configurations in $ \L $, $  \pr_\L\colon\bo\to\bo_\L:=\bo\cap \L\times Q $ for the projection from $ \bO $ to $ \bO_\L $ (similarly for unmarked configurations), $\boldsymbol{\Fcal}_\L^\prime=\boldsymbol{\Fcal}|_{\bO_\L} $ for the trace $\sigma$-algebra of $ \boldsymbol{\Fcal} $ on $ \bO_\L $, and  $ \boldsymbol{\Fcal}_\L=\pr_\L^{-1}\boldsymbol{\Fcal}_\L^\prime\subset\boldsymbol{\Fcal} $ for the $\sigma$-algebra of all events that happen in $ \L $ only. The reference measure on $ (\bO_\L,\boldsymbol{\Fcal}_\L^\prime) $ is $ \boldsymbol{\Pi}_\L^z:=\boldsymbol{\Pi}^z\circ\pr_\L^{-1} $. In a similar way we define the corresponding objects for unmarked configurations, $ \Pi^z, \Pi^z_\L, \O_\L, \pr_\L, \Fcal_\L^\prime $, and $ \Fcal_\L $. Finally, let $ \theta=(\theta_x)_{x\in\R^2} $ be the shift group, where $ \theta_x\colon\bO\to\bO $ is the translation of the  spatial component by the vector $ -x\in\R^2 $. Note that by definition, $ N_\Delta(\theta_x\bo)=N_{\Delta+x}(\bo) $ for all $ \Delta\subset\R^2 $. \medskip

Our Delaunay Potts model depends on the local geometry of the Delaunay triangulation. The Delaunay triangulation is a hypergraph structure comprising vertices, edges, and triangles and is defined as follows. The set $\Del(\omega)$ of Delaunay hyperedges of a given configuration $ \omega\in \Omega $ consists of all subsets $ \tau\subset\om $ for which there exists an open ball $ \Bb(\tau) \subset \R^2 $ with $ \partial \Bb(\tau)\cap\om=\tau $ which contains no points of $ \om $. For $ m=1,2,3  $, and $ \omega\in\Omega $, we write $$\Del_m(\omega)=\{\tau\in\Del(\omega)\colon \#\tau=m \} $$ for the set of Delaunay simplices with $ m $ vertices. It is possible that a Delaunay hyperedge $ \tau $ consists of four or more points on a sphere with no points inside. In fact, for this not to happen, we must consider configurations in general   position as in \cite{M94}. More precisely, this means that no four points lie on the boundary of a circle, no point lies inside a circumcircle for a triangle,  and every half-plane contains at least one point. Fortunately, this occurs with probability one for our Poisson reference measure, and in fact, for any stationary point process. Note that the open ball $ \Bb(\tau) $ is only uniquely determined when $ \#\tau=3 $ and $ \tau $ is affinely independent. In this case, $ \partial \Bb(\tau) $ is called the \emph{circumcircle} of $ \tau $ and its radius is called the \emph{circumradius} of $ \tau $ and is denoted $ \delta(\tau) $. Henceforth, for each configuration $ \om $, the associated set $\Del_3(\om)$ is known as the Delaunay triangulation of $\om$. It is uniquely determined and defines a triangulation of the convex hull of $ \om $ whenever $ \om $ is in general position (\cite{M94}). In a similar way one can define the marked Delaunay hyperedges $ \boldsymbol{\Del}_m(\bo) $,   where the Delaunay property refers to the spatial component only.

\medskip

In order to define our model in terms of a Gibbs distribution we introduce local versions of the defined Delaunay configurations. Given a configuration $ \om \in\Omega $ (or $ \bo$) we write $ \O_{\L,\om}=\{\zeta\in\O\colon \zeta\setminus\L=\om_{\L^{\rm c}}\} $ (resp. $ \bO_{\L,\bo} $) for the set of configurations which equal $ \om $ off $ \L $.  
The unmarked set
$$ 
\Del_{3,\L}(\zeta):=\{\tau\in\Del_3(\zeta)\colon\,\overline{\Bb(\tau)}\cap\L\not=\emptyset\} 
$$ 
and the corresponding marked set 
$$
\bDel_{3,\L}(\bzeta):=\{\btau\in\bDel_3(\bzeta)\colon \,\overline{\Bb(\tau)}\cap\L\not=\emptyset\} $$ contain the triangles in the Delaunay triangulation whose circumcircle intersects $ \L $. 

These are the triangles which can be removed from the triangulation by changing the configuration inside $\L$, that is to say
$$
\btau\in\bDel_{3,\L}(\bzeta)\Leftrightarrow \btau\in\bDel_3(\bzeta)\;\mbox{ and }\; \exists\,\bzeta^\prime\in\bO_\L\, \mbox{ s.t. } \btau\notin\bDel_{3}(\bzeta_{\L^{\rm c}}\cup\bzeta^\prime)\,.
$$

To see this, notice that if $\btau \in \bDel_{3,\Lambda}(\bzeta)$ and $\bx \in (\Lambda \cap \overline{\Bb(\tau)}) \times Q$ then $\btau \notin \bDel_3(\bzeta \cup \bx).$ On the other hand, suppose $\btau \in \bDel_3(\bzeta) \setminus \bDel_3(\bzeta_{\Lambda^c} \cup \bzeta')$ for some $\bzeta' \in \bO_\Lambda.$ If $\bzeta' = \emptyset$ then $\tau \cap \zeta_\Lambda \neq \emptyset$, else $\exists \bx \in \bzeta'$ such that $x \in \overline{\Bb(\tau)}$. In either case we must have $\overline{\Bb(\tau)} \cap \Lambda \neq \emptyset$, and so $\btau \in \bDel_{3,\Lambda}(\bzeta).$

\begin{figure}
    \centering
	\includegraphics[width=.5\linewidth]{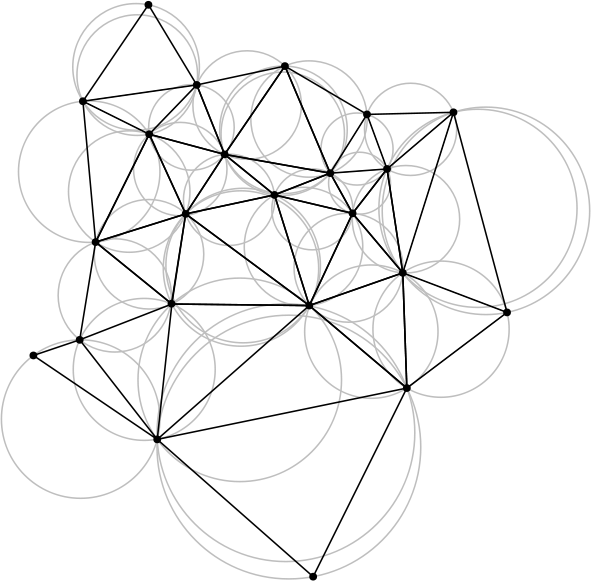}
	\caption{A Delaunay triangulation. The grey circles represent
	    the boundaries of the balls $B(\tau)$.}
    \label{}	
\end{figure}

For any triangle $ \tau\in\Del_3(\om) $ we denote its  area by $ \Asf(\tau) $ and its smallest interior angle by $ \alpha(\tau) $. The interaction is given by the following Hamiltonian in $ \L $ with boundary condition $ \bo\in\bO $,  
\begin{equation}\label{Hamiltonian}
H_{\L,\bo}(\bzeta):=\sum_{\tau\in\Del_{3,\L}(\zeta)}\;\Psi(\tau)+\sum_{\btau\in\bDel_{3,\L}(\bzeta)}\Phi_\beta(\Asf(\tau))(1-\delta_\sigma(\btau)),\quad\bzeta\in\O_{\L,\bo},
\end{equation}
where $ \Psi $ is the so-called \emph{background potential}
\begin{equation}
\Psi(\tau)=\begin{cases} 0 & \mbox{ if } \delta(\tau)\in (r,R)\;\mbox{ and }\; \alpha(\tau)>\alpha_0\,,\\
\infty & \mbox{ otherwise}\,,\end{cases} \quad 0<r<R<\infty, \ \alpha_0\in(0,\pi/3)\,.
\end{equation}
 Here  $ \Phi_\beta $ is the \emph{ferromagnetic type potential}  defined as a measurable function of the area $ \Asf(\tau) $ of a triangle defined for any $\beta\ge 0$,
\begin{equation}
\Phi_\beta(\Asf(\tau))=\log\Big(\frac{\Asf(\tau)+\beta}{\Asf(\tau)}\Big)\,,
\end{equation}
and 
$$ 
\delta_\sigma(\btau)=\begin{cases} 1 &, \mbox{ if } \sigma_{\btau}(x)=\sigma_{\btau}(y) =\sigma_{\btau}(z)\mbox{ for } \tau=\{x,y,z\},\\
0&, \mbox{ otherwise}.
\end{cases}
$$ 

\noindent Note the following scaling relation for the potential
\begin{equation}\label{scaleinvariance}
\Phi_\beta(L\Asf)=\Phi_{\beta/L}(\Asf),\qquad\mbox{for any } L>0\,,
\end{equation}
and that $ \Phi_\beta(\Asf(\tau))\to 0 $ when $ \Asf(\tau)\to\infty $ and $ \Phi_\beta(\Asf(\tau))\to \infty $ when $ \Asf(\tau)\to 0 $.
\smallskip

 Following \cite{DDG12} we define  the  partition function as
  $$
  Z_\L(\bo)=\int_{\bO_{\L,\bo}}\,{\rm e}^{-H_{\L,\bo}(\bzeta)}\,\Pi^z(\d\zeta^{\ssup{1}})\cdots\Pi^z(\d\zeta^{\ssup{q}}).
  $$
 The Gibbs distribution for the \emph{background} potential $ \Psi $ and the \emph{type} potential $ \Phi_\beta $, and $ z> 0 $ in $\L $ with boundary condition $ \bo $ is defined as
 \begin{equation}\label{Gibbsdist}
 \gamma_{\L,\bo}(A)=\frac{1}{Z_\L(\bo)}\int_{\bO_{\L,\bo}}\,\1_A(\bzeta\cup\bo){\rm e}^{-H_{\L,\bo}(\bzeta)}\,\Pi_\L^z(\d\bzeta),\quad A\in\boldsymbol{\Fcal}.
 \end{equation}
  It is evident from \eqref{Gibbsdist} that, for fixed $ \zeta\in\O_\L $, the conditional distribution of the marks of $ \bzeta= (\zeta^{\ssup{1}},\ldots,\zeta^{\ssup{q}})$  relative to $ \gamma_{\L,\bo} $ is that of a discrete Potts model on $ \zeta $ embedded in the Delaunay triangulation with position-dependent interaction between the marks. This justifies calling our model  \textit{Delaunay Potts model or Delaunay Widom-Rowlinson model}.
  
  \begin{definition}
  A probability measure $ \mu $ on $ \bO $ is called a Gibbs measure for the Delaunay  Potts model with activity $ z>0 $ and interaction type potential $ \Phi_\beta $  if 
 \begin{equation}\label{DLR}
   \E_\mu[f]=\int_{\bO}\,\frac{1}{Z_\L(\bo)}\int_{\bO_{\L,\bo}}f(\bzeta\cup\bo){\rm e}^{-H_{\L,\bo}(\bzeta)}\,\boldsymbol{\Pi}^z_\L(\d\bzeta)\mu(\d\bo)
 \end{equation} for every $ \L\Subset\R^2 $ and every measurable function $ f$.
  \end{definition}
  The equations in \eqref{DLR} are the DLR equations (after Dobrushin, Lanford, and Ruelle). They ensure that the Gibbs distribution in \eqref{Gibbsdist} is a version of the conditional probability $ \mu(A|\boldsymbol{\Fcal}_{\L^{\rm c}})(\bo) $. The measurability of all objects is established in \cite{E14,DDG12}.

\section{The Delaunay continuum random cluster measure}\label{sec-DelCluster}
We introduce the Delaunay continuum random cluster measure as an example of a hyperedge percolation model.  Hyperedge percolation models are created by taking random unmarked point configurations and declaring hyperedges of the associated hypergraph to be either `open' or `closed' according to some hyperedge process. A \emph{hypergraph structure} is a measurable subset $\Hcal $ of $ \O_{\rm f}\times\O  $ such that $ \tau\subset \omega $ for all $ (\tau,\omega)\in\Hcal $,  where $ \O_{\rm f} $ is the set of finite point configuration in $ \R^2 $. For a configuration $ \omega\in\O $ the pair $ (\omega,\Hcal)(\omega)) $ is called a \emph{hypergraph} where $ \omega $ is the set of \emph{vertices} and $ \Hcal(\omega):=\{\tau\in\O_{\rm f}\colon (\tau,\omega)\in\Hcal\} $ is the set of \emph{hyperedges}. In a similar way one defines the marked hypergraph structure $ \bHcal\subset\bO_{\rm f}\times\bO$. A notion of connectedness is obtained by declaring that two points are connected if one can travel between them via open hyperedges. We limit ourselves to unmarked hypergraphs for which each hyperedge contains the same number of points, that is there exists $ k\in\N $ such that $ \abs{\tau}=k $ for all $ (\tau,\omega)\in\Hcal $.  A configuration in this kind of model is a pair $ (\omega,\sigma_{\Hcal(\omega)}) $ where $ \omega\in\O $ and $ \sigma_{\Hcal(\omega)}\in\{0,1\}^{\Hcal(\omega)} $, with $ \sigma_{\Hcal(\omega)}(\tau)=1 $ signifying that $ \tau $ is open and $ \sigma_{\Hcal(\omega)}(\tau)=0 $ signifying that $ \tau $ is closed. Alternatively, we can represent the configuration as a pair $ G=(\omega,E) $ where $\omega\in\O $ and $E$ is locally finite set of hyperedges (the open ones). More precisely, $E\in\Escr $, where
$$ \Escr =\{E\subset E_{\R^2}\colon E\mbox{ locally finite}\}$$
and 
 $$
 E_{\R^2}=\{\tau\subset\R^2\colon \abs{\tau}=k\in\{1,2,3\}\}\,.
 $$ 
 
 The sample space is therefore $ \Gscr=\Omega\times\Escr $. In this formulation with $ G=(\omega,E)\in\Gscr $ the set $E$ represents the set of hyperedges which are considered to be open. Note that $ \Gscr $ also contains elements which do not belong to $ \Hcal $ since not every $ (\omega,E)\in\Gscr $ satisfies $ \tau\in E \Rightarrow \tau\subset\omega $.
 We equip $ \Escr $ with the $\sigma$-algebra $\Sigma $ generated by the counting variables $ N_\L\colon E\mapsto\abs{E_\L} $ for $ \L\Subset\R^{dk} $, similar to how we defined the $\sigma$-algebra $ \Fcal $ on $ \O$. As before, $ \Sigma $ is the Borel $\sigma$-algebra for the Polish topology on $ \Escr $ and thus the $\sigma$-algebra $ \Acal:=\Fcal\otimes\Sigma $ turns $ \Gscr $ into a Polish space.

 Let $ G=(\omega,E)\in\Gscr $. Two points $x,y\in\omega$ are \emph{adjacent} if there exists $ \tau\in E $ such that $ x,y\in\tau$. A \emph{path} connecting $x$ and $y$ is a sequence of points $ x_1,\ldots, x_n, x_i\in \omega, i=1,\ldots, n,n\in\N, $ with $ x_1=x $ and $x_n=y $ such that for all $ i\in\{1,\ldots, n-1\} $ there exists $ \tau_i\in E $ such that $ x_i,x_{i+1}\in \tau_i $. We say $x\in\omega$ and $y\in\omega $ belong to the same connected component of $G$ if there is a path connecting them.

 For the proof of Theorem~\ref{THM-main} we need the following version of the hyperedge model. Namely, we consider configuration with open and closed triangles (tiles).
 For $ \L\Subset\R^2 $  and parameters $ z $ and $ \Phi_\beta  $ we define a joint distribution of the Delaunay Potts model and an  process which we call \emph{Delaunay continuum random-cluster model}. The basic idea is to introduce random triangles  between points in the plane.   Let 
  $$
 E_{\R^2}=\{\tau\subset\R^2\colon \abs{\tau}=3\}\,.
 $$  
 be the set of all possible triangles  of points in $ \R^2 $, likewise, let $E_\L $ be  the set of all  triangle in $ \L $ and $ E_\zeta $ for the  set of triangles in $ \zeta\in\Omega_{\L,\om} $. We identify $ \om $ with $\om^{\ssup{1}} $ and $ \bo=(\om^{\ssup{1}},\emptyset,\ldots,\emptyset) $.  This allows only monochromatic boundary conditions whereas the general version involves the so-called Edwards-Sokal coupling (see \cite{GHM} for lattice Potts models). We restrict ourself to the former case for ease of notation.

 \noindent The joint distribution is built from the following two components.

 \noindent The \textit{point distribution}  is  given by the unmarked Gibbs distribution in \eqref{Gibbsdist}, i.e., the reference process is the Poisson process  $\Pi^{zq} $ for any  boundary condition $ \om\in\O $ and activity $ zq$ and 
 \begin{equation}\label{Gibbsdistun}
 \gamma_{\L,\omega}(A)=\frac{1}{Z_\L(\om)}\int_{\O_{\L,\om}}\,\1_A(\zeta\cup\om){\rm e}^{-H_{\L,\om}(\zeta)}\,\Pi_\L^{zq}(\d\zeta),\quad A\in\Fcal\,,
 \end{equation}  
 where
 $$
H_{\L,\om}(\zeta):=\sum_{\tau\in\Del_{3,\L}(\zeta)}\;\Psi(\tau)\,,\quad\zeta\in\O_{\L,\om}\,.
 $$ 
  
 \medskip

 \noindent The \textit{hyperedge drawing mechanism} respectively \emph{triangle drawing mechanism}. Given a point configuration $ \zeta\in\Omega_{\L,\om} $,   we let $ \mu_{\zeta,\L} $ be the distribution of the random hyperedge configuration $ \{\tau\in \Del_3(\zeta)\colon\upsilon(\tau)=1\} \in\Escr $ with the hyperedgeedge configuration $ \upsilon\in \{0,1\}^{\Del_3(\zeta)} $ having
 probability   
 $$
\prod_{\eta\in\Del_3(\zeta)}p(\tau)^{\upsilon(\tau)}(1-p(\tau))^{1-\upsilon(\tau)}
$$ with
\begin{equation}\label{tiledrawing}
p(\tau):=\P(\upsilon(\tau)=1)=\begin{cases} 1-\exp(-\Phi_\beta(\tau)) & \mbox{ if } \tau\in \Del_{3,\L}(\zeta)\,,\\
1 & \mbox{ if } \tau \in\Del_3(\zeta)\setminus\Del_{3,\L}(\zeta)\,.
\end{cases} 
\end{equation}
 
The measure $ \mu_{\zeta,\L} $ is a point process on $E_{\R^2} $ and is denoted the \emph{hyperedge drawing mechanism}. Note that $ \zeta\rightarrow\mu_{\zeta,\L} $ is a  probability kernel (see \cite{E14,AE16}). Let $ N_{\rm cc}(\om,E) $ denote the number of connected components in the hypergraph $ (\om,E) $. If
\begin{equation}
\Zsf_\L(\om):=\int\int\; q^{N_{\rm cc}(\zeta\cup \om_{\L^{\rm c}}, E)}\,\mu_{\zeta\cup\om_{\L^{\rm c}},\L}(\d E)\,\ex^{-H_{\L,\om}(\zeta)}\,\Pi_\L^{z}(\d\zeta)\in (0,\infty)
\end{equation}
then $ Z_{\L}(\om)\in (0,\infty) $ also, since $ q^{N_{\rm cc}(\om_{\L^{\rm c}}\cup \zeta, E)}\in (1,\infty) $. 

\begin{definition}[\textbf{Delaunay continuum random cluster measure}]
If $ \Zsf_\L(\om)\in (0,\infty) $ the \emph{Delaunay continuum random cluster measure} in $ \L\Subset\R^2 $ for $ \Psi,\Phi_\beta, z $ and boundary condition $ \om $ is the probability measure $ \Csf_{\L,\om}\in\Mcal_1(\Gscr,\Acal) $ defined by
\begin{equation}\label{def-drcm}
\Csf_{\L,\om}(A):=\frac{Z_\L(\om)}{\Zsf_\L(\om)}\;\1_A(\zeta\cup\om_{\L^{\rm c}},E)\,q^{N_{\rm cc}(\zeta\cup\om_{\L^{\rm c}},E)}\,\mu_{\zeta\cup\om_{\L^{\rm c}},\L}(\d E)\, \gamma_{\L,\om}(\d\zeta)\,,\quad A\in\Acal\,.
\end{equation}
\end{definition}

The final step is to obtain a representation measure for our Delaunay Potts model using the Delaunay continuum random cluster measure. We therefore need the final third component, the \emph{type} or \emph{mark picking mechanism}. 
 \noindent The \textit{type picking mechanism} for a given configuration $ \zeta\in\Omega_{\L,\om} $ is the distribution $ \lambda_{\zeta,\L} $ of the mark vector $\sigma_\zeta\in Q^\zeta$.  Here $ (\sigma_\zeta(x))_{x\in\zeta} $ are independent and uniformly distributed random variables on $ Q$ with $ \sigma_\zeta(x)=1 $ for all $ x\in\zeta_{\L^{\rm c}}=\om $.  The latter condition ensures that all points outside of $ \L$ carry the given fixed mark. We say that $ \lambda_{\zeta,\L} $ has a \emph{monochromatic} boundary condition since all points outside of $ \L $ have the same mark.
 
 Let $ B\in\bFcal\times\Sigma $ denote the event that any two points can only belong to the same connected component if they have the same mark, that is,
 $$
 B=\{(\bzeta,E)\in\bO\times\Escr\colon\sum_{\tau\in E}(1-\delta_{\sigma}(\btau))=0\}\,.
 $$ 
 If $ \om\in\O^*_\L:=\{\om\in\O\colon Z_\L(\om)\in (0,\infty)\} $ we define $ m_{\L,\om} $ on $ (\bO\times \Escr,\bFcal\times\Sigma) $ as the measure
 $$
 m_{\L,\om}(A):=\int\;\1_A(\bzeta,E) \mu_{\zeta\cup\om_{\L^{\rm c}},\L}(\d E)\,\lambda_{\zeta\cup\om_{\L^{\rm c}}}(\d\sigma_\zeta)\,\gamma_{\L,\om}(\d\zeta)\,,
 $$
 where $ \bzeta\equiv (\zeta,\sigma_\zeta) $.
 \medskip

\begin{definition}[\textbf{Delaunay random cluster representation measure}]

Delaunay random cluster representation measure in $ \L\Subset\R^2 $ for $ \Psi, \Phi_\beta,z $, and boundary condition $ \om\in\O^*_\L $ is the probability measure
on $ (\bO\times\Escr,\bFcal\times \Sigma) $ defined by
$$
\Psf_{\L,\om}(\cdot):=m_{\L,\om}(\cdot |B)\,.
$$

\end{definition}

The random cluster representation measure is a joint construction of the Delaunay continuum Potts distribution and the Delaunay continuum random cluster distribution. The former can be obtained if one only looks at the particle positions and their types and disregards the hyperedges. Alternatively, the latter can be obtained by are ignoring the type of each particle. These statements are formalised in the next proposition. They are very similar to Propositions \cite[2.1]{GH96} and \cite[2.2]{GH96}. The proofs in our context do not differ in any notable way so they are omitted.

Let $ \rho_1 $ and $ \rho_2 $ denote the projections from $ \bO\times\Escr $ to $ \bO $ and $ \O\times\Escr $ respectively.

\begin{proposition}
Let $ \L\Subset\R^2 , z>0 ,\om\in\O^*_\L $ and $ \bo\equiv(\om,\sigma_\om) $ where $ \sigma_\om\equiv 1 $. Then the following holds
\begin{enumerate}[(a)]
\item $\Psf_{\L,\om}\circ \rho_1^{-1}=\gamma_{\L,\om}\,,$

\medskip

\item $ \Psf_{\L,\om}\circ\rho_2^{-1}=\Csf_{\L,\om}\,.$
\end{enumerate}
\end{proposition}

 \section{Gibbs Tessellations with geometric hardcore conditions and type interaction}\label{sec-existence}

  \begin{theorem}[\textbf{Existence of Gibbs tessellations}]\label{Propexist}
 For $ \beta>0, \alpha_0\in (0,\pi/3) $ and given $ 0<r<R<\infty $ there exists a $ z_0=z_0(\beta,R,r,\alpha_0) $ such that  
 for any $z>z_0 $ there exists at least one translation-invariant Gibbs measure for the Delaunay Potts)  model with potentials $ \Psi $ and $ \Phi_\beta $. 
 
  \end{theorem}

 \begin{remark}[\textbf{Gibbs measures}]
 The proof is using the so-called pseudo-periodic configurations (see Appendix~\ref{pseudoperiodic} or \cite{DDG12}) and properties of the potential $ \Phi_\beta $. Existence of Gibbs measures for related different  Delaunay models have been obtained in \cite{BBD99, Der08,DG09}. Note that for $q=1$ our models have no marks and Gibbs measures do exist as well (\cite{DDG12,Der08}).  In that case only background potential $\Psi $ is relevant which is equivalent to the hardcore potential $ \Psi^\prime $ introduced in \cite{Der08}. Namely,
 $$
 \Psi^\prime(\tau,\eta)=\begin{cases} 0 & \mbox{ if } \delta(\tau)<R  \; \mbox{ or }\; \ell(\eta)> r\,,\\ \infty & \mbox{ otherwise}\,,\end{cases}
 $$  
 where $ \eta\in\Del_2 $ and $ \tau\in\Del_2 $. Here, $\ell(\eta) $ denotes the length of the Delaunay edge $ \eta\in\Del_2 $, i.e., $ \ell(\eta)=\abs{x-y} $ if $ \eta=\{x,y\} $.
 
  \hfill $ \diamond$
 \end{remark}

To prove that a translation-invariant Delaunay continuum Potts measure exists (Proposition \ref{Propexist}) we will show that the appropriate adapted versions of the conditions outlined in \cite{DDG12} are satisfied. These conditions are the \textit{range condition} (R), \textit{stability} (S), and \textit{upper regularity} (U). (R) is satisfied for the same reason as was given for the Delaunay Potts models discussed in \cite{AE16}, and (S) is trivially satisfied since the potentials $\Psi$ and $\Phi_\beta$ are positive. Upper regularity is composed of three separate conditions: \textit{uniform confinement, uniform summability,} and \textit{strong non-rigidity}. To verify these conditions we need to introduce the notion of pseudo-periodic configurations. The definition we use here suits our purposes but is less general than the definition given in \cite{DDG12}. Restrictions will be placed on the parameters $\ell$ and $\rho$ throughout the course of the proof.
\begin{definition}{(Pseudo-periodic configurations.)}	
	\label{marked pseudo-periodic full def}
	Let $\rho,\ell>0$ and
	\begin{align*}
	    M = \small 
	    \begin{pmatrix}
		    \ell & \ell/2 \\
		    0 & \sqrt{3}\ell/2
	    \end{pmatrix}.
    \end{align*} 
    Define for $k = (i,j) \in \Z^2$ the cells
	\begin{align}
		\label{define cells}
		\Delta_k = \Delta_{i,j} := \left\{ Mx \in \R^d : x - k \in [-1/2, 1/2)^2 \right\},
	\end{align}
	which form a partition of $\R^2$ into rhombi, and the set
	\begin{align*}
	    \bGamma:= \{ \bo \in \bO_{\Delta_{0,0}} : \omega = \{ x \} \text{ for some } x \in B(0,\rho \ell) \},
	\end{align*}
	where $B(0,\rho \ell)$ is the ball of radius $\rho \ell$ around the origin.
	The configurations that belong to the set
	\begin{align}
		\label{pseudo-periodic def 1}
		\overline{\bGamma} := \{ \bo \in \bO : \theta_{Mk} (\bo_{\Delta_k}) \in \bGamma \text{ for all } k \in \Z^2 \}.
	\end{align}
	are called \textit{pseudo-periodic}.
\end{definition}

In our case, to show uniform confinement it is sufficient to find some $k>0$ such that for any $\Delta \Subset \R^2, \bo \in \overline{\bGamma}, \bzeta \in \bO_\Delta, \btau \in \bDel_{3,\Delta}(\bzeta \cup \bo_{\Delta^c})$ and $\bo' \in \bO$ satisfying $\bo'_{\Delta \oplus k} = \bo_{\Delta \oplus k}$ we have $\btau \in \bDel_3(\bzeta \cup \bo_{\Delta^c}').$ Clearly if $\btau \in \bDel_{3,\Delta}(\bzeta \cup \bo_{\Delta^c})$ and $\Bb(\tau) \subset \Delta \oplus k$ then $\btau \in \bDel_3(\bzeta \cup \bo_{\Delta^c}')$, so it is sufficient to show that there exists some $k$ such that $\Bb(\tau) \subset \Delta \oplus k$ for any $\btau \in \bDel_{3,\Delta}(\bzeta \cup \bo_{\Delta^c})$. Such a value $k$ must exist since $\bo$ contains a point in each cell $\Delta_{i,j}$. 

The potentials $\Phi_\beta$ and $\Psi$ are positive, so the requirement of uniform summability is that
\begin{align}
    \label{uniform summability}
    c_{\bGamma} := \sup_{\bo \in \overline{\bGamma}} \sum_{{\substack{\btau \in \bDel_3(\bo): \\ \tau \cap \Delta_{0,0} \neq \emptyset}}} \frac{\Psi(\tau) + \Phi_\beta(\Asf(\tau))(1 - \delta_\sigma(\btau))}{|\tilde{\btau}|} < \infty,
\end{align}
where $\tilde{\btau} := \{ k \in \Z^2 : \tau \cap \Delta(k) \neq \emptyset \}$. We first verify that if $\bo \in \overline{\bGamma}$ and $\btau \in \bDel_3(\bo)$ then $\Psi(\tau)=0$. 

\begin{figure}
	\centering
	\includegraphics[width=0.8\linewidth]{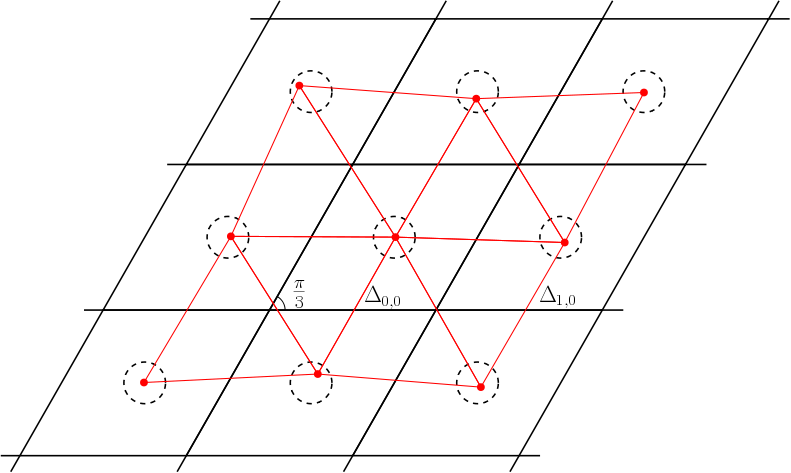}
	\caption{A pseudo-periodic configuration with small $\rho$.}
	\label{tessellation_2}	
\end{figure}

For small enough $\rho$ (it suffices that $\rho < 1/(2\sqrt{3})$), each point $\bx \in \bo \in \overline{\bGamma}$ has 6 neighbours and the Delaunay triangulation of $\bo$ becomes a perturbed triangular lattice (see Figure \ref{tessellation_2}.)
In this case the length of each edge lies in the interval $(\ell(1 - 2\rho), \ell(1 + 2\rho))$. Thus by the law of cosines 
\begin{align*}
	\cos (\alpha(\tau)) &\leq \frac{2\ell^2(1+2\rho)^2 - \ell^2(1-2\rho)^2}{2\ell^2(1-2\rho)^2} \\
	&= \left( \frac{1 + 2\rho}{1 - 2\rho} \right)^2 - \frac{1}{2}
\end{align*}
for all $\tau \in \Del_3(\omega).$ The roots of the quadratic $(1 + 2\rho)^2 - (\frac{1}{2} + \cos(\alpha_0))(1 - 2\rho)^2$ are $\frac{(1 \pm (\frac{1}{2} + \cos(\alpha_0))^{\frac{1}{2}})^2}{2\cos(\alpha_0)-1}$, and so since $\cos(\alpha_0) > \frac{1}{2}$ we have the following result.

\begin{lemma}
    \label{hardcore angle satisfied for pseudo periodic}
    If $\rho < \frac{(1 - (\frac{1}{2} + \cos(\alpha_0))^{\frac{1}{2}})^2}{2\cos(\alpha_0)-1} \wedge \frac{1}{2\sqrt{3}}$ then for any pseudo-periodic configuration $\bo \in \overline{\bGamma}$,
    \begin{align*}
        \alpha(\tau) > \alpha_0 \text{ \normalfont for all } \tau \in \Del_3(\omega).
    \end{align*}
\end{lemma}

The perimeter of a triangle $\tau$ is bounded above by $3\sqrt{3}\delta(\tau)$ (with equality when $\tau$ is equilateral), so if $\btau \in \bDel_3(\bo)$ and $\bo \in \overline{\bGamma}$ then
\begin{align*}
    \frac{\ell (1-2\rho)}{\sqrt{3}} \leq \delta(\tau)
\end{align*}
On the other hand, the circumradius of a triangle with area $\Asf$ and edge lengths $a,b$ and $c$ is $\frac{abc}{4\Asf}$, so we have 
\begin{align*}
    \frac{(1-2\rho)}{\sqrt{3}} \leq \frac{\delta(\tau)}{\ell} \leq \frac{\ell^2(1 + 2\rho)^3}{4\Asf(\tau)}.
\end{align*}
Now the lower bound 
\begin{align}
    \label{area lower bound}
    \Asf(\tau)
	&\geq \sqrt{3\ell(1-2\rho)(2\ell(1-2\rho)-\ell(1+2\rho))^3} \nonumber \\
	&= \ell^2 \sqrt{3(1-2\rho)(1-6\rho)^3}
\end{align}
can be obtained using Heron's formula, so if we further assume that $\rho < \frac{1}{6},$ we have
\begin{align}
	\label{upper bound circumradius pseudo-periodic}
	0 < L(\rho) := \frac{1-6\rho}{\sqrt{3}} \leq \frac{\delta(\tau)}{\ell} \leq \frac{(1 + 2\rho)^3}{\sqrt{3(1-2\rho)(1-6\rho)^3}} \leq \frac{(1 + 6\rho)^3}{\sqrt{3}(1-6\rho)^2} =: U(\rho) .
\end{align}
These inequalities are used to prove the following result.
\begin{proposition}
    \label{proposition - both hardcore satisfied}
   Let
   $$
   \rho_0(r,R,\alpha_0) := \frac{R^{1/3} - r^{1/3} }{6 \left( R^{1/3} + r^{1/3} \right)} \wedge \frac{(1 - (\frac{1}{2} + \cos(\alpha_0))^{\frac{1}{2}})^2}{2\cos(\alpha_0)-1}\,.
 $$  If $\rho < \rho_0(r,R,\alpha_0) $  then $\frac{r}{L(\rho)} < \frac{R}{U(\rho)}$. Furthermore, if $\ell \in \left( \frac{r}{L(\rho)}, \frac{R}{U(\rho)} \right)$ then all pseudo-periodic configurations $\bo \in \overline{\bGamma}$ satisfy
    \begin{align}
        \label{both hardcore satisfied}
        \delta(\tau) \in (r,R) \text{\normalfont \ and } \alpha(\tau) > \alpha_0 \text{ \normalfont for all } \tau \in \Del_3(\omega),
    \end{align}
    and therefore $\Psi(\tau)=0$.
\end{proposition}
\begin{proofsect}{Proof}
    The first part is just a simple rearrangement:
    \begin{align*}
        & \rho < \frac{R^{1/3} - r^{1/3} }{6 \left( R^{1/3} + r^{1/3} \right) }, \\
        & \implies (1+6\rho) r ^{1/3} < (1 - 6\rho)R^{1/3}, \\
        & \implies \frac{r}{R} < \left( \frac{1-6\rho}{1+6\rho} \right)^3  = \frac{L(\rho)}{U(\rho)}, \\
        & \implies \frac{r}{L(\rho)} < \frac{R}{U(\rho)}.
    \end{align*}
    Since $\rho < \frac{R^{1/3} - r^{1/3} }{6 \left( R^{1/3} + r^{1/3} \right)} < \frac{1}{6} < \frac{1}{2\sqrt{3}}$ we can apply Lemma \ref{hardcore angle satisfied for pseudo periodic} and inequality (\ref{upper bound circumradius pseudo-periodic}) to obtain  
    \begin{align*}
        r < \ell L(\rho) \leq \delta(\tau) \leq \ell U(\rho) < R
    \end{align*}
    and $\alpha(\tau) > \alpha_0$ for all $\tau \in \Del_3(\omega)$ when  $\ell \in \left( \frac{r}{L(\rho)}, \frac{R}{U(\rho)} \right)$.
    \qed

\end{proofsect}

\medskip

If the conditions of Proposition \ref{proposition - both hardcore satisfied} are satisfied then 
\begin{align*}
	 c_{\bGamma} \leq \sup_{\bo \in \overline{\bGamma}} \sum_{{\substack{\btau \in \bDel_3(\bo): \\ \tau \cap \Delta_{0,0} \neq \emptyset}}} \frac{\Phi_\beta(\Asf(\tau)))}{3}.
\end{align*} 
Thus applying inequality (\ref{area lower bound}) yields
\begin{align}
    \label{c_gamma upper bound}
	c_{\bGamma} \leq 2 \log \left( 1 + \frac{\beta}{\ell^2 \sqrt{3(1-2\rho)(1-6\rho)^3}} \right) < \infty,
\end{align} 
which concludes the proof of (\ref{uniform summability}).

Finally, the requirement of strong non-rigidity is that $$\ex^{z|\Delta_{0,0}|} \bPi^z_{\Delta_{0,0}}(\bGamma) > \ex^{c_{\bGamma}}\,. $$  Since $$\ex^{z|\Delta_{0,0}|} \bPi^z_{\Delta_{0,0}}(\bGamma) = z|B(0,\rho \ell)|\,, $$  this inequality is satisfied when
\begin{align*}
    z >  \frac{1}{\pi \rho^2 \ell^2} \left( 1 + \frac{\beta}{\ell^2 \sqrt{3(1-2\rho)(1-6\rho)^3}} \right) ^2.
\end{align*}
Thus far we have established ranges for each parameter for which (R), (S) and (U) are satisfied, and therefore a translation-invariant Delaunay continuum Potts measure exists. The precise statement is given below.

\begin{proposition}
    \label{existence theorem hardcore}
    If $\beta>0$,
    \begin{align*}
        \rho < \rho_0(R,r,\alpha_0) :=  \frac{R^{1/3} - r^{1/3} }{6 \left( R^{1/3} + r^{1/3} \right) } \wedge \frac{(1 - (\frac{1}{2} + \cos(\alpha_0))^{\frac{1}{2}})^2}{2\cos(\alpha_0)-1},
    \end{align*}
    \begin{align*}
        \ell \in \left( \frac{r}{L(\rho)}, \frac{R}{U(\rho)} \right),
    \end{align*} 
    and 
    \begin{align}
        \label{z0 prime definition}
        z > z'_0(\beta, \rho, \ell) := \frac{1}{\pi \rho^2 \ell^2} \left( 1 + \frac{\beta}{\ell^2 \sqrt{3(1-2\rho)(1-6\rho)^3}} \right)^2,
    \end{align} 
    then there exists a translation-invariant Delaunay continuum Potts measure for $\bDel_3,z,\beta,\Phi_\beta$ and $\Psi$.
\end{proposition}
Proposition \ref{Propexist} is a corollary of Proposition \ref{existence theorem hardcore}.

\section{Phase Transitions for sufficiently large activities and large potential parameter}\label{sec-phase}

  A \textit{phase transition}  is said to occur if there exists more than one Gibbs measure for the Delaunay Potts model. The following  theorem shows that this happens for sufficiently large  activities $ z $ and large potential parameter $ \beta >0 $.   Note that $ \beta $ is a parameter for the type interaction and not the usual inverse temperature.

  \begin{theorem}[\textbf{Phase transition}]\label{THM-main}
 Let $ \alpha_0\in(0,\sin^{-1}(3/64) $, and $ 0<r<3/64 R <\infty $. There exists $\beta_0=\beta_0(q,R,r,\alpha_0) >0 $ and $ z_0=z_0(\beta,q,R,r,\alpha_0)>0 $ such that for all $ \beta>\beta_0 $ and  $z>z_0 $   there exist at least $q$ different  translation-invariant Gibbs measures for the Delaunay continuum Potts model. 
  \end{theorem}

\begin{remark}
\begin{enumerate}[(a)]

\item Theorem~\ref{THM-main} actually establishes a break of the symmetry in the type distribution.

\smallskip

\item Theorem~\ref{THM-main} also holds for any type potential  $\Phi_\beta^{\ssup{\gamma}}  $ of the form 
$$
\Phi_\beta^{\ssup{\gamma}}(\Asf(\tau)):=\log\Big(1+\beta\Asf(\tau)^{-\gamma}\Big),\qquad \gamma>0\,,
$$
with different lower bounds $ z_0 $ and $ \beta_0 $.

\end{enumerate}   \hfill $ \diamond $
\end{remark}

\begin{remark}[\textbf{Free energy and Uniqueness of Gibbs measures}]
\begin{enumerate}[(a)]
\item One may wonder if the phase transition manifest itself thermodynamically by a non-differentiability ("discontinuity") of the free energy (pressure).  We refer the interested reader to \cite{AE16} for more details and references.

\smallskip

 \item  To establish uniqueness of the Gibbs measure in our Delaunay Potts model one can use the Delaunay random-cluster measure $ \Csf_{\L_n,\omega} $ defined in \eqref{def-drcm} in Section~\ref{sec-DelCluster}. This would require to prove absence of percolation in the model given by   $ \Csf_{\L_n,\omega} $ and in future work we address this question. \end{enumerate}
\hfill $ \diamond$ 
 \end{remark}

Before turning to details of the proof of Theorem~\ref{THM-main} we outlined first why percolation in the Delaunay continuum random cluster measure introduced in Section~\ref{sec-DelCluster} leads to the dominance of a given type and hence the breaking of the uniform type distribution.  The main body of the proof is then devoted to establish this percolation by coarse-graining methods. 

For $ \L\Subset\R^2 $ being some union of cells $ \Delta_i\subset \L$ to be specified later, let $ N_{\L,i}(\bo) $ denote the number of particles in $ \L $ carrying the mark $ i\in Q $ and let $ N_{\Delta\leftrightarrow\L^{\rm c}}(\om,E) $ denote the number of particles in $ \Delta$ which are connected to $ \L^{\rm c} $:
$$
\begin{aligned}
N_{\L,i}(\bo)&:=\abs{\{x\in\om_\L\colon \sigma_\om(x)=i\}}\,,\\
N_{\Delta\leftrightarrow\L^{\rm c}}(\om,E) &:=\abs{\{x\in \om_\Delta\colon \,\exists\, \mbox{ a path in }\, (\om,E) \,\mbox{ from } x\in\Delta \mbox{ to some } y\in\L^{\rm c}\}}\,.\end{aligned}
$$

The following proposition shows how the dominance of a given type is linked to percolation of the Delaunay continuum random cluster measure. This goes back to \cite{GH96} and we provide and adapted version following \cite{AE16,AE19}.
 
 \begin{proposition}[\textbf{\cite{AE16,AE19}}]\label{Propkey}
Let  $ \L\Subset\R^2 $ be some union of cells $ \Delta_i\subset \L$, $\om\in\O_\L^* $  and choose $\bo $ monochromatic such that $ \sigma_{\bo}(x)=1 $ for all $ x\in\om $. Then
$$
\int\;\Big(qN_{\L,1}(\bzeta\cup\bo_{\L^{\rm c}})-N_\L(\bzeta\cup\bo_{\L^{\rm c}})\Big)\,\gamma_{\L,\bo}(\d\bzeta)=(q-1)\int\;  N_{\Delta\leftrightarrow\L^{\rm c}}(\zeta\cup \om_{\L^{\rm c}},E) \Csf_{\L,\om}(\d\zeta,\d E)\,.
$$
 
 \end{proposition}

 The key task in the whole proof is to establish percolation for the Delaunay continuum random cluster measure, that is, the obtain a strict lower bound of the right hand side uniformly for all possible unions of cells. 
  
  \begin{proposition}\label{Prop-per}
  Suppose all the assumptions hold and that $ z $ and $ \beta $ are sufficiently large. Suppose that $ \L $ is a finite union of cells $ \Delta_{k,l} $ defined in \eqref{cell} Appendix A. Then there exists $\eps>0 $ such that
  \begin{equation}\label{pereq}
  \int\, N_{\Delta\leftrightarrow\L^{\rm c}}(\zeta\cup\om_{\L^{\rm c}}, E) \,\Csf_{\L,\om}(\d\zeta,\d E) \ge \eps
  \end{equation} for any cell $\Delta=\Delta_{k,l} $, any finite union $ \L $ of cells and any pseudo-periodic boundary condition $ \om \in\O_\L^*$ such that $ \bo\in\bGamma$.  
  \end{proposition}

Before proving the pivotal statement in Proposition~\ref{Prop-per} we sketch how to complete the proof of the breaking of uniform type distribution in Theorem~\ref{THM-main}.

\begin{proofsect}{Proof of Theorem~\ref{THM-main}}
We give only a sketch of the ideas going back to \cite{GH96} and in our context from \cite{AE16,AE19}. In the following 
$$ 
\L\equiv\L_n=\bigcup_{(k,l)\in\{-n,\ldots,n\}^2}\;\Delta_{k,l}\,.
$$
 We shall construct a sequence of probability measures $ (P_N^{\ssup{1}})_{n\in\N} $ on $ (\bO,\bFcal) $ such that the following holds. 
 \begin{enumerate}[(i)]
 \item $P_n^{\ssup{1}} $ are invariant under the skewed lattice transformations $ (\theta_x)_{x\in M\Z^2} $\,,
 
 \item For any type $ i\in Q $ and $ \Delta\subset \L_n $ and boundary condition $ \om\in\O^*_{\L_n} $ with $ \bo\in\bGamma $ and $ \sigma_{\bo}(x)=1 $ for all $x\in \om $,
 $$
 \int\:\Big(qN_{\Delta,i}-N_\Delta\Big)\,\d P_n^{\ssup{1}}=\int\;\Big(qN_{\Delta,i}-N_\Delta\Big)\,\d\gamma_{\L,\om}\,,
 $$
   
 \item   $ (P_N^{\ssup{1}})_{n\in\N} $ has a subsequence which converges locally\footnote{In this instance the local convergence topology is the weak* topology generated by the set of local and tame real-valued functions on $\bO $. These are the functions $f$ which ate $ \bFcal_\Delta$-measurable and satisfy $\abs{f(\bo)}\le a \abs{\om_{\Delta}}+b $ for some $ \Delta\Subset\R^2 $ and $ a,b\in\R $\,, see Appendix~\ref{AppC}.} to some measure $ P^{\ssup{1}} $, and after spatially averaging the measure $ P^{\ssup{1}}(\cdot|\{\emptyset\}^{\rm c}) $, one obtains a translation-invariant Gibbs measure $ \widetilde{P}^{\ssup{1}} $. 
 
  \end{enumerate} 
 
 The local convergence in conjunction with the uniform bound in Proposition~\ref{Prop-per} and Proposition~\ref{Propkey} then implies that for all $ \Delta\Subset\R^2 $,
 $$
  \int\:\Big(qN_{\Delta,i}-N_\Delta\Big)\,\d \widetilde{P}^{\ssup{1}}\ge (q-1)\eps>0 \,, $$
 and since the Gibbs distribution $ \gamma_{\L,\bo} $ is invariant under permutations of the remaining type $ j\in Q\setminus\{1\} $, we have
 $$
 \int\; N_{\Delta,1} \, \d \widetilde{P}^{\ssup{1}}>\int\; N_{\Delta,j} \d \widetilde{P}^{\ssup{1}}\,,\quad \mbox{ for all } j\in Q\setminus\{1\}\,.
 $$
 This shows that the symmetry of the type distribution is broken. Finally, in the same way, for each $j\in Q\setminus\{1\} $ it is possible to obtain a translation-invariant Gibbs measure $  \widetilde{P}^{\ssup{j}}$ in which the preferred mark is $j $ instead of $1$. This concludes the proof.

\qed
\end{proofsect}

   \subsection{Hyperedge percolation to site percolation}\label{sec-siteversion}
We shall show how the uniform lower bound condition in Proposition~\ref{Prop-per} can be achieved. We will construct a continuum site percolation model $ \Csf^{\ssup{\text{site}}}_{\L,\omega}$ in which the event that $\Delta$ is connected to $\L^{\rm c}$ is smaller that it is with respect to the Delaunay continuum random cluster measure $\Csf_{\L,\omega}.$ We  then use a coarse graining argument to bound this event from below. The new percolation model shares the same particle distribution as $\Csf_{\L,\omega} $ define din Section~\ref{sec-DelCluster}. 

First we must introduce the notion of stochastic dominance between probability measures. A function $f \colon \Ecal \rightarrow \R$ is said to be \textit{increasing} if $f(A) \leq f(B)$ whenever $A \subset B$. For two probability measures $\mu_1,\mu_2$ on $(\Ecal,\Sigma),$ we say that $\mu_1$ \textit{stochastically dominates} $\mu_2$ and write $\mu_1 \succcurlyeq \mu_2$ if $\mu_1(f) \geq \mu_2(f)$ for all increasing functions $f \colon \Ecal \rightarrow \R$.

The measure $   \Csf^{\ssup{\text{site}}}_{\L,\omega} $ is  be defined as a measure on $(\bO,\bFcal)$ where $q=2$, although instead of the mark space $\{1,2\}$ we will use the mark space $\{0,1\}.$ Points with mark $1$ are considered to be `open' and points with mark $0$ are `closed.' A \textit{path} in $\bo \in \bO$ connecting $\bm{x}$ and $\bm{y}$ is a sequence of points $(\bm{x}_i)_{i=1}^n \subset \bo$ with $n \in \N, \bm{x}_1 = \bm{x}$ and $\bm{x}_n = \bm{y}$ such that $\sigma_\omega(x_i) = 1$ for all $i \in [n]$ and there exists $\btau_j \in \bDel_3(\bo)$ such that $\bm{x}_j,\bm{x}_{j+1} \in \btau_j$ for all $j \in [n-1]$. The event that $\Delta$ is connected to $\L^{\rm c}$ is the following:
\begin{align*}
    \{\Delta \leftrightarrow \L^{\rm c} \} := 
    \left\{ 
        \bo \in \bO \ \bigg|
        \begin{array}{cc}
            \text{There exists a path } (\bm{x}_i)_{i=1}^n  \\
            \text{ in } \bo \text{ with } x_1 \in \omega_\Delta \text{ and } x_n \in \omega_{\L^{\rm c}}.
        \end{array}
    \right\}
\end{align*}

Let $M_{\L,\omega}$ denote the marginal distribution $\Csf_{\L, \omega}(\cdot,\Ecal)$. We can then write 
\begin{align*}
	\Csf_{\L,\omega}(d \omega', \d E) = \mu^q_{\omega',\L}(\d E) M_{\L,\omega} (\d \omega'). 
\end{align*}
where
\begin{align}
	\label{tile drawing with q}
	\mu^q_{\omega',\L}(\d E) := \frac{q^{K(\omega',E)} \mu_{\omega',\L}(\d E)}{\int q^{K(\omega',E)} \mu_{\omega',\L}(\d E)}.
\end{align}

Let $\widehat{\Hcal}$ denote an unmarked hypergraph structure which is a subset of $\Del_3$, and pick $\hat{p} \in [0,1].$ 

For $\omega \in \O$, let $\widehat{\mu}_\omega$ denote the distribution of the random hyperedge configuration $\left\{ \tau \in \Del_3(\omega): \xi_\tau = 1 \right\} \in \Ecal,$ where $(\xi_\tau)_{\tau \in \widehat{\Hcal}(\omega)}$ are independent Bernoulli random variables such that $\xi_\tau = 1$ with probability $\hat{p} \1_{\widehat{\Hcal}(\omega)}(\tau).$ In other words, each hyperedge $\tau \in \Del_3(\omega)$ is declared open independently with probability $\hat{p} \1_{\widehat{\Hcal}(\omega)}(\tau),$ and closed otherwise. 

In addition, define $$\widehat{\omega} := \{ x \in \omega : \exists \tau \in \widehat{\Hcal}(\omega) \text{ with } x \in \tau\} $$ and let $\widehat{\lambda}_\omega$ denote the distribution of the mark vector $\sigma_\omega \in \{0,1\}^\omega$ where $(\sigma_\omega(x))_{x \in \omega}$ are independent and identically distributed such that $\sigma_\omega(x) = 1$ with probability $p \1_{\widehat{\omega}}(x)$.
We can now define our site percolation measure on $(\bO,\bFcal)$ to be
\begin{align*}
	 \Csf^{\ssup{\text{site}}}_{\L,\omega}(\d \omega', \d \sigma_\omega) :=  \widehat{\lambda}_\omega(\d \sigma_{\omega'}) M_{\L,\omega} (\d \omega').
\end{align*}
The result regarding connectivity is the following. For the full proof see \cite[Proposition 2.18]{E14} or \cite{AE16, AE19}.
\begin{proposition}[\textbf{\cite[Proposition~2.18]{E14},\cite{AE16,AE19}}]
	\label{site comparison}
	If $\mu^{q}_{\omega,\L} \succcurlyeq \hat{\mu}_\omega$ then for all $\Delta \subset \L \Subset \R^2$,
	\begin{equation}
	 	\int N_{\Delta \leftrightarrow \L^{\rm c}} \d \Csf_{\L,\omega} \geq  \Csf^{\ssup{\text{site}}}_{\L,\omega}(\Delta\leftrightarrow\L^{\rm c})\,.
\end{equation}
		\end{proposition}

Let $\widehat{\mu}_{\omega,\Delta}$ denote the measure for which each hyperedge $\tau \in \Del_3(\omega)$ is opened independently with probability 
\begin{align*}
    \widehat{p}(\tau)\1_{\widehat{\Hcal}(\omega) \cap \Del_{3,\L}}(\omega)(\tau) + \1_{\Del_3(\omega) \setminus \Del_{3,\L}(\omega)}(\tau),
\end{align*} where $\widehat{p}: \Ecal_{\R^2,3} \rightarrow [0,1].$ With respect to both $\mu^{q}_{\omega,\L}$ and $\widehat{\mu}_{\omega,\L}$, the status of all but finitely many hyperedges are fixed. In this case, one can show that $\mu^{q}_{\omega,\L} \succcurlyeq \widehat{\mu}_{\omega,\L}$ if for all $\tau \in \Del_3(\omega)$, the \textit{comparison inequalities}
\begin{align}
	\label{comparison inequalities}
	\frac{p_\Delta(\tau)}{q^2(1-p_\Delta(\tau))} \geq \frac{\hat{p}(\tau)}{(1-\hat{p}(\tau))} 
\end{align}
are satisfied (with the convention that $\frac{p}{1-p} = \infty$ when $p=1$) by applying the same method as in the case where the hypergraph is finite. For the proof in the case of a finite hypergraph, see \cite[Proposition 2.3]{E14} (which generalises a result originally proven in \cite{F72}). Employing a coupling argument similar to that used in Lemma \ref{percolation greater than bernoulli} we see that $\widehat{\mu}_{\omega,\L} \succcurlyeq \widehat{\mu}_\omega$. Therefore, to show that $\mu^{q}_{\omega,\L} \succcurlyeq \widehat{\mu}_\omega$ we need only verify (\ref{comparison inequalities}).

\noindent\textbf{Proof of Proposition~\ref{Prop-per}}: We split the proof in several steps and Lemmata below.  To achieve the uniform lower bound on the right hand side of \eqref{pereq} in Proposition~\ref{Prop-per} we show this lower bound for the random cluster measure $ \Csf^{\ssup{\rm site}}_{\L,\om} $ defined in Section~\ref{sec-siteversion}. As this percolation measure is stochastically smaller we obtain the desired uniform bound ones we have established the bound for the new site process $ \Csf^{\ssup{\rm site}}_{\L,\om} $. The pivotal idea here is the employ a \emph{coarse graining} procedure, that is, in Section~\ref{coarse graining section} we will pave the plane $ \R^2 $ with different cell systems ranging between three different scales depending on the scale given by the background potential. The coarse graining procedure allows to define good cells and ultimate allows to prove the percolation for the measure  $ \Csf^{\ssup{\rm site}}_{\L,\om} $ by comparing to site percolation on $ \Z^2 $ where each good cell with centre in $ \Z^2 $ is declared open.  This is done in Section~\ref{sec-percolation} resulting in the percolation for $ \Csf^{\ssup{\rm site}}_{\L,\om} $ for certain values of the parameter involved. This then leads to percolation with uniform bounds for the Delaunay continuum random cluster measure $ \Csf_{\L,\om} $. The main technical work is to find and then to gauge  the correct coarse graining procedure adapted to our problem. This requires control of the particle distribution once an additional particle is inserted into a configuration. The challenge here is that due to the \emph{non-hereditary} feature of the Delaunay tessellations that an additional particle can lead to an invalid configuration. The details for this delicate control of the particle distribution is in Section~\ref{sec-augmentation}.

\qed

\subsection{Augmenting a configuration by a particle}\label{sec-augmentation}
The situation here differs from the case of a classical many-body interaction since adding a point $x$ does not merely result in additional interaction terms representing the interaction between $x$ and the other particles. Instead, when the particle $x$ is added, some hyperedges are created and others are destroyed. The result of this is that for $\omega' \in \O_\Delta$ and $x \in \Delta$, both $H_{\L,\omega}(\omega')$ and   $H_{\L,\omega}(\omega' \cup \{x\})$ may contain terms which are not present in the other. It is thus possible (depending on the background potential $\Psi$) that 
$$
H_{\L,\omega}(\omega') = \infty \; \mbox{ and }\; H_{\L,\omega}(\omega' \cup \{x\}) < \infty\,. 
$$ 
In this case the function $\ex^{-H_{\L,\omega}(\cdot)}$ is called \textit{non-hereditary} \cite[Definition 10.4.IV]{DV08}. We start with the \textit{point insertion lemma} which expresses that the circumradius of each `new' triangle (when $x$ is added) is less than the circumradius of particular `old' triangles. The point insertion lemma formalises an argument found in section 12.2.6 (page 462) of \cite{Li99}, with some details filled in.

\begin{figure}[htbp]
	\centering
	\subfigure[$T_{x_0, \omega}^- = \Del_3(\omega) \setminus \Del_3(\omega \cup x_0)$]{\includegraphics[scale=.6]{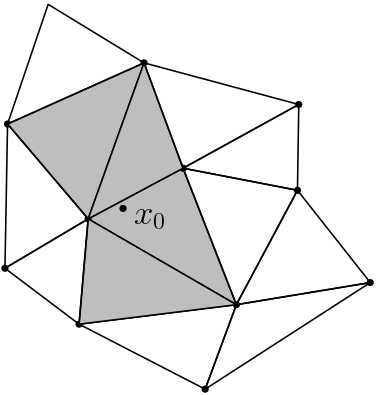}}

	\subfigure[$T_{x_0, \omega}^+ = \Del_3(\omega \cup  x_0 ) \setminus \Del_3(\omega)$ ]{\includegraphics[scale=0.6]{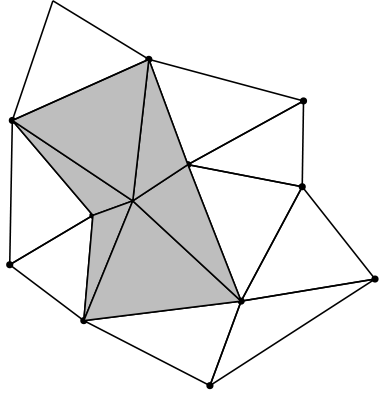}}
				
	\caption{Augmenting the configuration $\omega$ by a point $x_0$.}
	\label{adding x_0 picture}	
\end{figure}

The following sets contain the tiles that remain intact, the ones that are created, and those that are destroyed when adding the point $x_0$ to an unmarked configuration $\omega$. From now on we will write $\omega \cup x_0$ rather than $\omega \cup \{x_0\}.$
\begin{align*}
	T_{x_0, \omega}^{\text{ext}} &:= \Del_3(\omega) \cap \Del_3(\omega \cup x_0) = \{\tau \in \Del_3(\omega) : x_0 \notin \overline{B(\tau)}\}, \\
	T_{x_0, \omega}^+ &:= \Del_3(\omega \cup x_0 ) \setminus \Del_3(\omega) = \{\tau \in \Del_3(\omega \cup x_0) : x_0 \in \tau \}, \\
	T_{x_0, \omega}^- &:= \Del_3(\omega) \setminus \Del_3(\omega \cup x_0) = \{\tau \in \Del_3(\omega) : x_0 \in \overline{B(\tau)}\}.
\end{align*}
An example is shown in figure \ref{adding x_0 picture}. The area covered by the triangles in $T_{x_0, \omega}^-$ (or $T_{x_0, \omega}^+$) is shown in grey and referred to as the \textit{Delaunay cavity} created by $x_0$. 

\begin{lemma} \emph{\textbf{(Point insertion lemma.)}}
    \label{point insertion}
    Suppose $\tau = \{ x_0,y,z \} \in \Del_3(\omega \cup x_0)$ and let $\tau_1,\tau_2$ denote the two triangles in $\Del_3(\omega)$ which have $\{y,z\}$ as a subset. Then
    \begin{align*}
        \delta(\tau) \leq \max(\delta(\tau_1),\delta(\tau_2)).
    \end{align*}
\end{lemma}

\begin{figure}[htbp]
	    \centering
	    \subfigure[$\d (u,C_\tau)^2 =  \d (u,C_{\tau_2})^2 + \d (C_\tau, C_{\tau_2})^2 + 2\d (C_\tau, C_{\tau_2}) \d (C_{\tau_2}, p)$]{\includegraphics[scale=0.35]{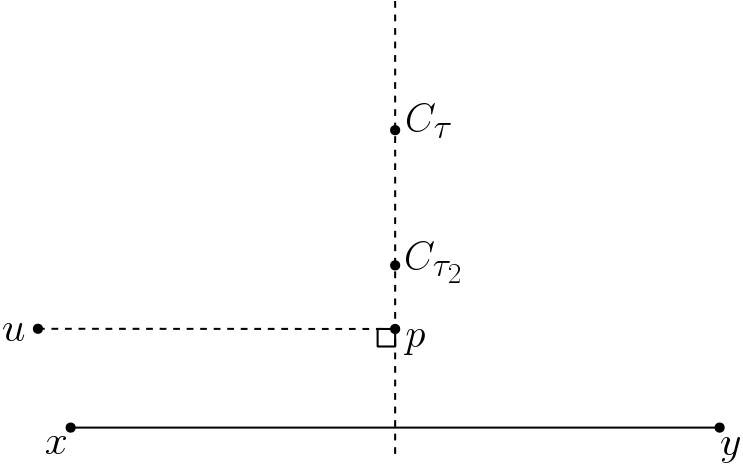}}
	       
	    \vspace{1cm}

        \subfigure[ $\d (u,C_\tau)^2 \leq \d (u,C_{\tau_2}) + \d (C_\tau, C_{\tau_2})^2.$]{\includegraphics[scale=0.4]{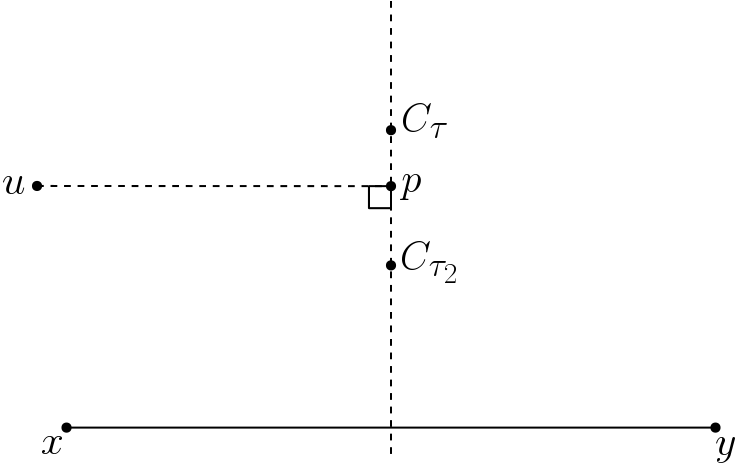}}

 \vspace{1cm}
 
  \subfigure[ $\d (u,C_\tau) \leq \d (u,C_{\tau_2}).$ ]{\includegraphics[scale=0.2]{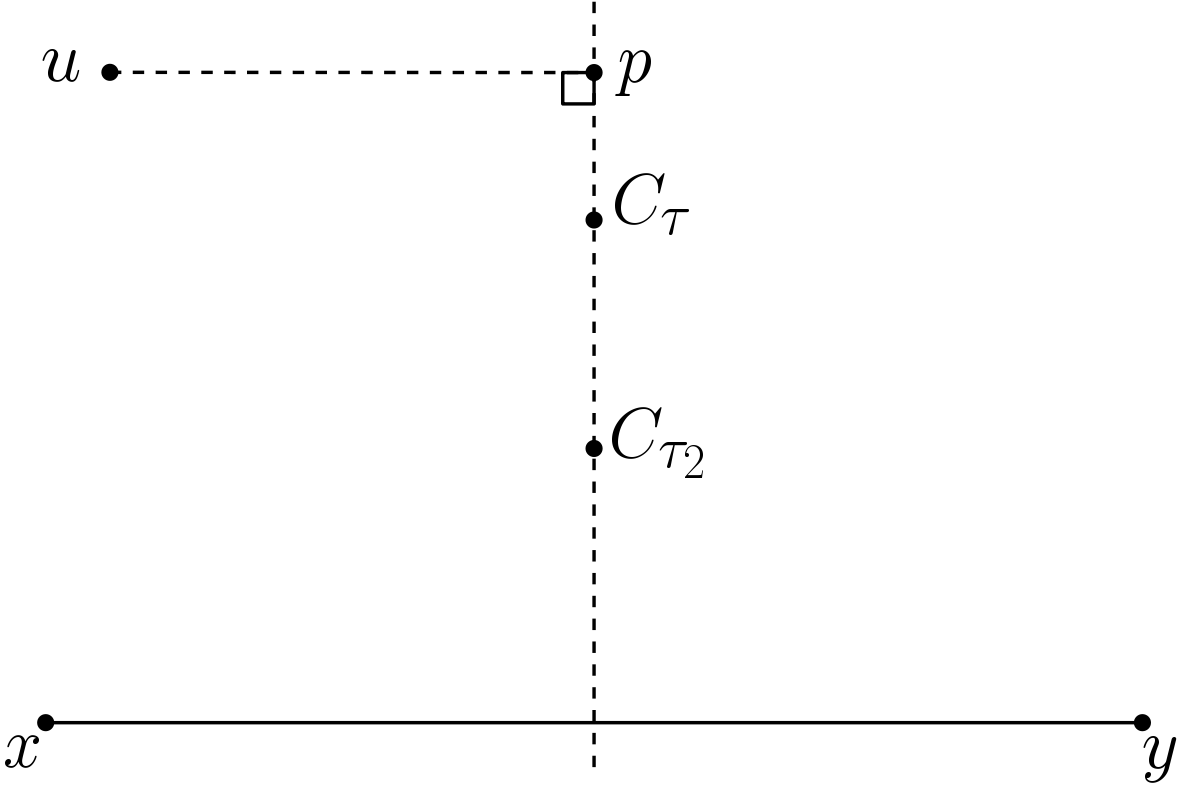}}
\caption{}
        \label{diagrams - distance between u and Ctau}	
    \end{figure}

\begin{proofsect}{Proof}
	Without loss of generality, let $\tau_1 = \{v,y,z\} \in T^-_{x_0, \omega}$ and $\tau_2 = \{u,y,z\} \in T^{\text{ext}}_{x_0, \omega}$ (see figure \ref{Delaunay cavity}). Let $C_\tau$ denote the circumcentre of $\tau$. Consider the two half planes separated by the line $\overleftrightarrow{xy}$ passing through $x$ and $y$. We will show that if $C_\tau$ is in the same half-plane as $v$ then $\delta(\tau) \leq \delta(\tau_1)$ and if $C_\tau$ is in the same half plane as $u$ then $\delta(\tau) \leq \delta(\tau_2)$. 
	\begin{figure}
        \centering
        \subfigure{	\includegraphics[scale=0.25]{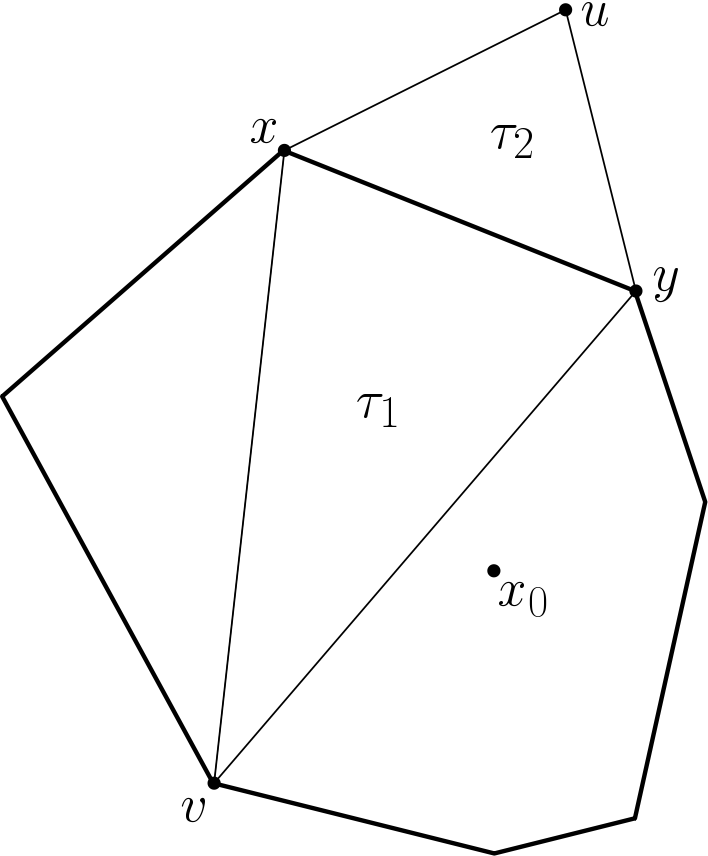}}
	      
 \subfigure{\includegraphics[scale=0.5]{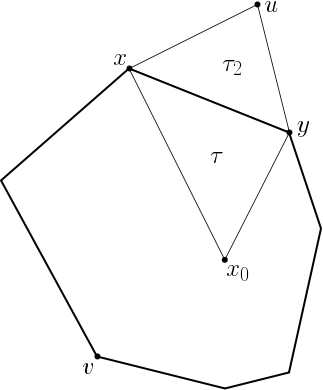}}
    	 \caption{}   
        \label{Delaunay cavity}	
    \end{figure}
	
	In the former case, the angle $\theta_1$ subtended by the chord $\overline{xy}$ at $v$ is less than the angle $\theta_2$ subtended by $\overline{xy}$ at $x_0$. This can be seen by extending the line $\overline{yx_0}$ until it intersects the circumcircle of $\tau_1$ (figure \ref{subtended angles}), which is possible since $x_0$ lies inside the circumcircle of $\tau_1.$

	\begin{figure}
    	\centering
    	\includegraphics[width=0.35\linewidth]{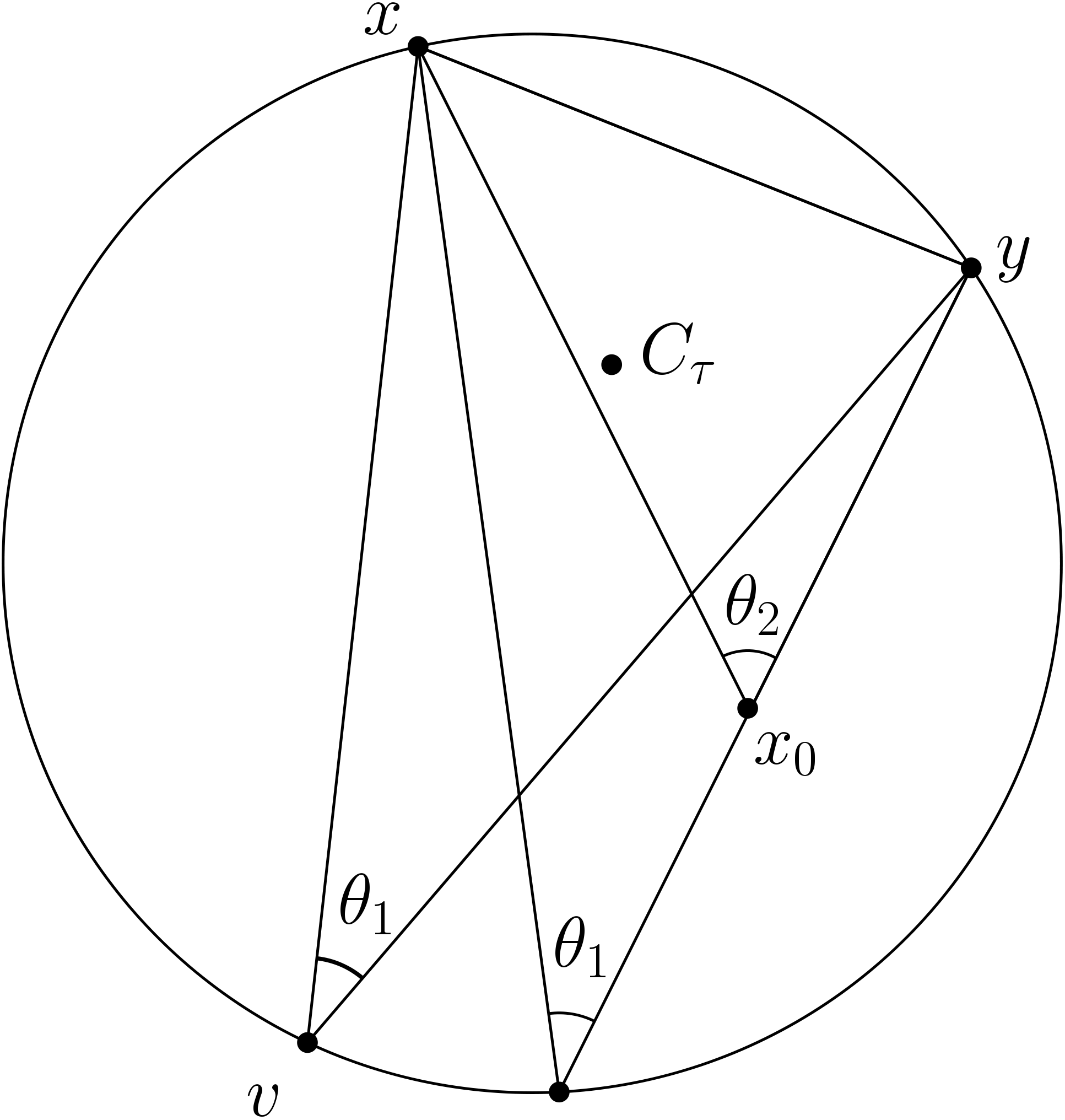}
        \caption{The angle $\theta_1$ subtended by the chord $\overline{xy}$ at $v$ is less than the angle $\theta_2$ subtended by $\overline{xy}$ at $x_0$.}
	    \label{subtended angles}
    \end{figure}

    \begin{figure}
        \centering
		\includegraphics[width=0.35\linewidth]{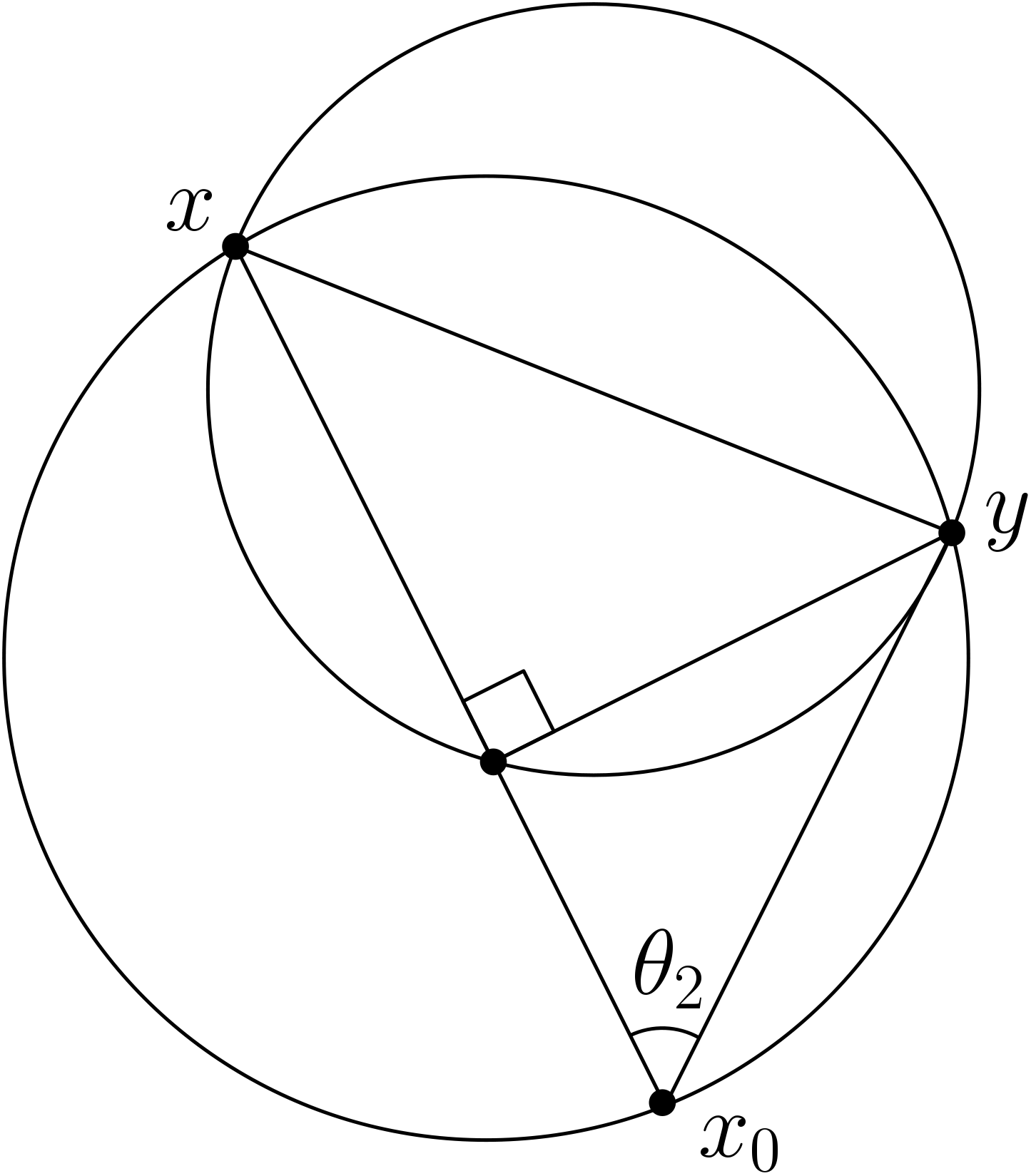}
    	\caption{There is a right angle subtended by $\overline{xy}$ at some point along $\overline{xx_0}$ or $\overline{yx_0}$.}
	    \label{subtended right angle}
    \end{figure}

	Since $C_\tau$ lies on the same side of $\overline{xy}$ as $v$, there is a right angle subtended by $\overline{xy}$ at some point along $\overline{xx_0}$ or $\overline{yx_0}$ (figure \ref{subtended right angle}), so $\theta_2 \leq \frac{\pi}{2}$. From the relationship $\theta_1 \leq \theta_2 \leq \frac{\pi}{2}$ we can conclude that
	\begin{align*}
		\delta(\tau) = \frac{\d (x,y)}{2 \sin(\theta_2)} \leq \frac{\d (x,y)}{2 \sin(\theta_1)} = \delta(\tau_1).
	\end{align*}
	    
	Now suppose that $C_\tau$ is in the same half plane as $u$. If $C_{\tau_2}$ is farther away from $\overline{xy}$ than $C_\tau$, then it is clear that $\delta(\tau) = \d (x,C_\tau) \leq \d (x,C_{\tau_2}) = \delta(\tau_2)$, as required. 
	
	Conversely, if $C_{\tau_2}$ is between $\overline{xy}$ and $C_{\tau}$ then $\delta(\tau_2) \leq \delta(\tau)$. In fact
	\begin{align}
	    \label{point insertion - circumradius inequality}
	    \delta(\tau_2) \leq \delta(\tau) < \d (u,C_\tau)
	\end{align} 
	since the point $u$ is outside of $\overline{B(\tau)}$. But we also have (see figure \ref{diagrams - distance between u and Ctau})
	\begin{align*}
	    \delta(\tau)^2 
	    &= \d(x,C_\tau)^2 \\
	    &= \d (x,C_{\tau_2})^2 + \d (C_\tau, C_{\tau_2})^2 + 2\d (C_\tau, C_{\tau_2}) \d \left( C_{\tau_2}, \frac{x + y}{2} \right)\\
	    &= \d (u,C_{\tau_2})^2 + \d (C_\tau, C_{\tau_2})^2 + 2\d ( C_\tau, C_{\tau_2}) \d \left( C_{\tau_2}, \frac{x + y}{2} \right)\\
	    &\geq \d(u,C_\tau)^2.
	\end{align*}
	which contradicts (\ref{point insertion - circumradius inequality}).
	\qed     
	\end{proofsect}

\noindent Let $\omega \in \O^*_\L$. Recall that $M_{\L,\omega}$ denotes the marginal distribution $\Csf_{\L,\omega}(\cdot,\Ecal).$ The Radon-Nikodym density of $M_{\L,\omega}$ with respect to $\gamma_{\L,\omega}$ is 
\begin{align*}
    h_{\L,\omega}(\omega') := \1_{ \{ \omega_{\L^{\rm c}} = \omega'_{\Delta^c} \} } \frac{Z_\L(\omega)}{\mathsf{Z}_\L(\omega)} \int q^{N_{\rm cc}(\omega',T)} \mu_{\omega',\Delta}(\d T).
\end{align*}
The next lemma gives a lower bound on the ratio
\begin{align*}
    \frac{h_{\L,\omega}(\omega' \cup x_0)}{h_{\L,\omega}(\omega')},
\end{align*}
which is known as the \textit{Papangelou conditional intensity}. Recall that according to the hyperedge drawing mechanism $\mu_{\omega,\L}$ each edge in $\Del_3(\omega)$ is opened independently according to the probabilities given in \eqref{tiledrawing}. Let $\mu^{\text{ext}}_{x_0,\omega},\mu^+_{x_o,\omega}$ and $\mu^-_{x_0,\omega}$ denote the measures which open the edges in $T^{\text{ext}}_{x_0,\omega}, T^+_{x_0,\omega}$ and $T^-_{x_0,\omega}$ respectively with the same probabilities. Then 
\begin{align*}
	\mu^{\text{ext}}_{x_0,\omega} \otimes \mu^+_{x_0,\omega} = \mu_{\omega \cup x_0,\L} \qquad \text{ and } \qquad \mu^{\text{ext}}_{x_0,\omega} \otimes \mu^-_{x_0,\omega} = \mu_{\omega,\L}.
\end{align*}

\begin{lemma}
	\label{hardcore_papangelou_intensity}
	Suppose that $\omega \in \O^*_\L$ with $\omega_{\L^{\rm c}} = \omega'_{\L^{\rm c}},$ and $x_0 \in \L \setminus \omega'$. If 
	$$
	H_{\L,\omega}(\omega'_\L), H_{\L,\omega}(\omega_\L' \cup x_0 ) < \infty\,,$$ then 
	\begin{align*}
		\frac{h_{\L,\omega}(\omega' \cup x_0)}{h_{\L,\omega}(\omega')} \geq q^{1-\frac{2\pi}{\alpha_0}}.
	\end{align*}
\end{lemma}
\begin{proofsect}{Proof} 
The structure of this proof is the same as \cite[Lemma 2.3]{AE19} but the details are slightly different since we are dealing with triangle interactions rather than edge interactions.
	\begin{align*}
		\frac{h_{\L,\omega}(\omega' \cup x_0)}{h_{\L,\omega}(\omega')} 
		&= \frac{\int q^{N_{\rm cc}(\omega' \cup x_0, T)} \mu_{\omega' \cup x_0, \Delta}(\d T)}{\int q^{N_{\rm cc}(\omega', T)} \,\mu_{\omega', \Delta}(\d T)} \\
		&= \frac{\int q^{N_{\rm cc}(\omega' \cup x_0, T_1 \cup T_2) - N_{\rm cc}(\omega', T_1)} q^{N_{\rm cc}(\omega', T_1)} \,\mu^{\text{ext}}_{x_0,\omega'}(\d T_1) \mu^+_{x_0,\omega'}(\d T_2)}{\int q^{N_{\rm cc}(\omega', T_3 \cup T_4) - N_{\rm cc}(\omega', T_3)} q^{N_{\rm cc}(\omega', T_3)} \,\mu^{\text{ext}}_{x_0,\omega'}(\d T_3) \mu^-_{x_0,\omega'}(\d T_4)}. 
	\end{align*}
	Opening more triangles can only reduce the number of connected components, so
	\begin{align*}
	    N_{\rm cc}(\omega', T_3 \cup T_4) \leq N_{\rm cc}(\omega', T_3).
	\end{align*}
	Furthermore, since $H_{\L,\omega}(\omega_\L' \cup x_0 ) < \infty$ the point $x_0$ is connected to at most $\frac{2\pi}{\alpha_0}$ other points in the graph $(\omega' \cup x_0, \Del_2(\omega' \cup x_0)),$ so 
	\begin{align*}
	    N_{\rm cc}(\omega' \cup x_0, T_1 \cup T_2) - N_{\rm cc}(\omega', T_1) \geq 1 - \frac{2\pi}{\alpha_0}.
	\end{align*}
	Therefore
	\begin{align*}
		\frac{h_{\L,\omega}(\omega' \cup x_0)}{h_{\L,\omega}(\omega')} 
		&\geq \frac{\int q^{1-\frac{2\pi}{\alpha_0}} q^{N_{\rm cc}(\omega', T_1)} \mu^{\text{ext}}_{x_0,\omega'}(\d T_1) \mu^+_{x_0,\omega'}(\d T_2)}{\int q^{N_{\rm cc}(\omega', T_3)} \mu^{\text{ext}}_{x_0,\omega'}(\d T_3) \mu^-_{x_0,\omega'}(\d T_4)} \\
		&= q^{1-\frac{2\pi}{\alpha_0}}. \qedhere
	\end{align*}
\qed

\end{proofsect}

\subsection{Coarse graining}\label{coarse graining section}
In order to prove the existence of a uniform lower bound on $\Csf_{\L_n,\omega}^{\ssup{\text{site}}}(C(k,m) \leftrightarrow \L_n^{\rm c})$, we will devise a criterion by which, according to the underlying configuration $\bo$, each cell $C(k,m)$ (defined in \ref{define cells}) is declared open or closed. This criterion will be devised in such a way that there exists an infinite connected component containing a point in $C(k,m)$ if $C(k,m)$ belongs to an infinite connected component of open boxes. We call this procedure of moving from points to cells \textit{coarse graining}. Formally, for each $n$ we will construct a map $X_n: \bO \rightarrow \{0,1\}^{\Z^2}$ where $\bo \in \{C(k,m) \leftrightarrow \L_n^c\}$ if $(k,m)$ belongs to an infinite open cluster in $X_n(\bo).$ The cell $C(k,m)$ is considered to be open if $X_n(\bo)(k,m) = 1.$ The desired lower bound will then be obtained by making a stochastic comparison between the law of $X_n$ and a Bernoulli product measure using Corollary ~\ref{mixed-site bond positive percolation probability}.
Let $M$ and $\bGamma$ be as in Definition~\ref{marked pseudo-periodic full def}, and let the parameters $\rho,\ell$ and $z$ satisfy the requirements of Proposition \ref{existence theorem hardcore} with $zq$ in place of $z$. Then \textbf{(U)} is satisfied (in addition to \textbf{(R)} and \textbf{(S)}). The same argument as in the previous section can be used to show that these conditions are also satisfied in the unmarked regime with respect to $\Psi$ instead of $\Phi_\beta$, and $\Gamma := \{ \omega \in \O_C : \omega = \{ x \} \text{ for some } x \in B \}$ instead of $\bGamma$. The situation is simpler since in this case $c_\Gamma = 0.$ As the boundary condition is pseudo-periodic we have $\bGamma \subset \bO^*_{\L_n}$ and $\Gamma \subset \O^*_{\L_n}$ for all $n \in \N$, where $\L_n = \bigcup_{k,m \in \{-n...,n\}}C(k,m)$.

The cells $(C(k,m))_{k,m \in \Z}$ form a partition of the plane into rhombi of length $\ell$. Let us split each cell into 64 smaller sub-cells of length $\ell/8$ denoted $(C_{k,m}^{i,j})_{0 \leq i,j \leq 7}$, where
\begin{align*}
    C_{k,m}^{i,j} := \left\{ Mx \in \R^2 : x - (k,m) \in \left[\frac{i-4}{8},\frac{i-3}{8} \right ) \times  \left[\frac{j-4}{8}, \frac{j-3}{8} \right) \right\}.
\end{align*}
Let $F_{k,m}$ denote the event that there is least one particle in each sub-cell of $C(k,m)$ and $O_{k,m}$ denote the event that additionally all points in $C(k,m)$ are open:
\begin{align*}
    F_{k,m} :=& \ \bigcap_{0 \leq i,j \leq 7} \left\{ \bo \in \O : | \omega \cap C_{k,m}^{i,j}| \geq 1 \right\}. \\
    O_{k,m} :=& \ \{ \bo \in \bO : \sigma_{\omega} (x) = 1 \text{ for all } x \in \omega_{C_{k,m}} \}.
\end{align*}
The map $X_n: \bO \rightarrow \{0,1\}^{\Z^2}$ is constructed by opening the sites $(k,m)$ inside $\{-n,...,n\}^2$ for which $F_{k,m} \cap O_{k,m}$ occurs, and opening the sites outside $\{-n,...,n\}^2$ for which $O_{k,m}$ occures, i.e 
\begin{align*}
    X_n(\bo)(k,m) := 
	\begin{cases} 
		\1_{F_{k,m} \cap O_{k,m}}(\bo) & \mbox{ if } |k|,|m| \leq n  \\ 
		\1_{O_{k,m}}(\bo) & \mbox{ otherwise.}
	\end{cases}	
\end{align*}
$X_n$ is therefore a stochastically decreasing sequence. To complete the definition of $\Csf_{\L_n,\omega}^{\ssup{\text{site}}} $, let $\widehat{\Hcal} = \Del_3$  and 
\begin{align*}
    \hat{p}(\tau) \equiv \hat{p} = \frac{1}{\frac{3\sqrt{3}}{4\beta} q^2 R^2 + 1}.
\end{align*} 
For a given circumradius, the triangle $\tau$ with maximal area is the equilateral triangle, for which $\Asf(\tau) = \frac{3\sqrt{3}}{4}\delta(\tau)$. Therefore, if $H_{\L_n,\omega}(\omega') < \infty$ then 
\begin{align*}
    \hat{p} &\leq \frac{1}{q^2\beta^{-1} \Asf(\tau) + 1}\\
    \implies \frac{\hat{p}}{1-\hat{p}} &\leq \frac{\beta}{q^2\Asf(\tau)}
\end{align*}
for all $\tau \in \Del_{3,\Delta}(\omega_{\L_n^c} \cup \omega').$ The comparison inequalities (\ref{comparison inequalities}) are satisfied since
\begin{align*}
    \frac{p_\L(\tau)}{q^2(1-p_\L(\tau))} = \frac{1 - \ex^{-\Phi_\beta(\tau)}}{q^2 \ex^{-\Phi_\beta(\tau)}} = \frac{\ex^{\Phi_\beta(\tau)}-1}{q^2} = \frac{\beta}{q^2\Asf(\tau)},
\end{align*}
for all $\tau \in \Del_{3,\L}(\omega_{\L_n^{\rm c}} \cup \omega')$, and therefore $\mu^{q}_{\omega',\L_n} \succcurlyeq \widehat{\mu}_{\omega'}$ almost surely with respect to $M_{\L,\omega}.$ Hence the premises of Theorem~\ref{THM-main} and Proposition~\ref{Prop-per} are satisfied, so it remains to show that there exists $c>0$ such that $  \Csf_{\L_n,\omega}^{\ssup{\text{site}}}(C(k,m) \leftrightarrow \L_n^c) \geq c$ for all $n$ and $\omega \in \Gamma.$ The following lemma shows that a uniform lower bound on the percolation probability of the law of $X_n$ is sufficient.

\begin{lemma}
    If $\omega \in \Gamma, \sigma_\omega(x) \equiv 1$, $H_{\L_n,\omega}(\omega') < \infty$ and $X_n(\bo_{{\L_n}^{\rm c}} \cup \bo') \in \{(k,m) \leftrightarrow \infty\}$, then $\bo_{{\L_n}^{\rm c}}' \cup \bo' \in \{C(k,m) \leftrightarrow \L_n^{\rm c} \}.$ Therefore
    \begin{align*}
       \Csf_{\L_n,\omega}^{\ssup{\rm site}} (C(k,m) \leftrightarrow \L_n^{\rm c}) \geq \mathcal{L}_{X_n}((k,m) \leftrightarrow \infty) 
    \end{align*}
    where $\mathcal{L}_{X_n}$ is the law of $X_n$.
\end{lemma}

\begin{figure}
    \centering
    \includegraphics[width=0.75\linewidth]{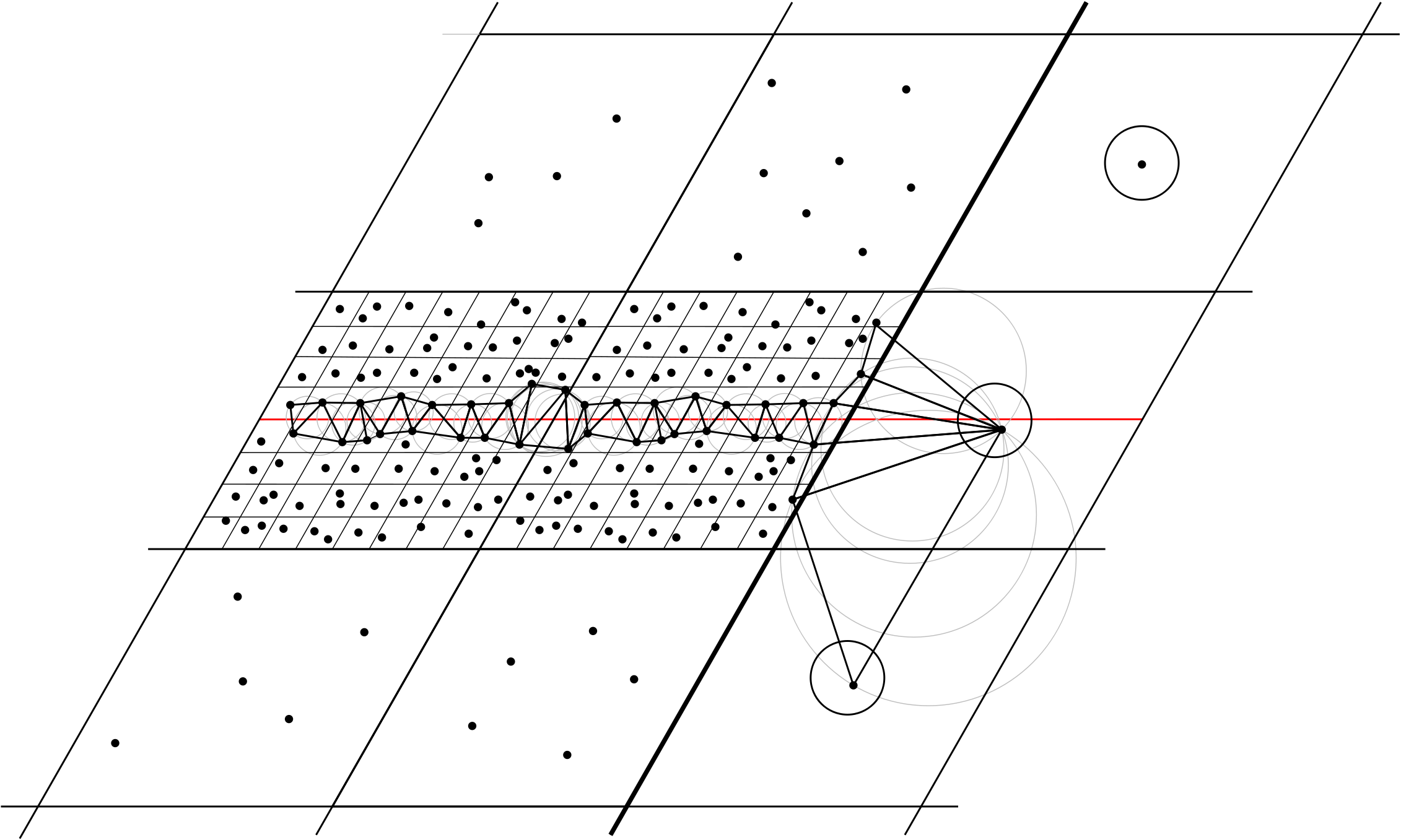}
	\caption{An illustration of two open cells meeting the boundary of $\L_n$, which is represented by the bold line. The open cells in $\L_n$ have at least one point in each of their 64 sub-cells.}
	\label{path between open cells}
\end{figure}

\begin{proofsect}{Proof}
    For $k \in \{-n,...,n-1\},|m| \leq n$ let $x_{k,m},x_{k+1,m} \in \omega_{\L^c} \cup \omega'$ denote the points whose Voronoi cells contain the centers of $C(k,m)$ and $C(k+1,m).$ The Voronoi cell associated to $x$ is given by 
    \begin{align*}
        \Vor_{\omega_{\L^c} \cup \omega'}(x) := \{ z \in \R^2 : |x-z| \leq |w - z| \text{ for all } w \in \omega_{\L^c} \cup \omega' \}.
    \end{align*}
    If $X_n(\bo_{\L_n^c} \cup \bo')(k,m) = X_n(\bo_{\L_n^c} \cup \bo')(k+1,m) = 1$ then $\bm{x}_{k,m}$ is connected to $\bm{x}_{k+1,m}$ via a path whose points are located in \begin{align*}
        \left\{ \bigcup_{2 \leq i \leq 7, 2 \leq j \leq 5} C_{k,m}^{i,j} \right\} \cup \left\{ \bigcup_{0 \leq i \leq 5, 2 \leq j \leq 5} C_{k+1,m}^{i,j} \right\}.
    \end{align*} 
    This can be seen via the same argument as \cite[Lemma 2.7, step (iv)]{AE19}. In fact, the points can be joined via a path whose points all have Voronoi cells intersecting the line segment between the centers of the cells $C(k,m)$ and $C(k+1,m).$  
    
    Furthermore, the same applies if $|m| \leq n$ and $X_n(\bo_{\L_n^c} \cup \bo')(n,m) = X_n(\bo_{\L_n^c} \cup \bo')(n+1,m) = 1$; there is a path between $\bm{x}_{n,m}$ and the point in $\bo_{C(n+1,m)}$ (recall that $\rho < \frac{1}{6}$) via a path whose points all have Voronoi cells intersecting the line segment between the centers of the cells $C(n,m)$ and $C(n+1,m).$ Figure \ref{path between open cells} shows an example of a path passing through two open cells and across the boundary of $\L_n.$ 
    
    Note that we have only discussed horizontal crossings between cells $C(k,m)$ and $C(k+1,m)$. The proof for vertical crossings can be performed similarly. It is now clear that if there is an infinite path $(k_r,m_r)_{r=1}^s$ in $X_n(\bo_{{\L_n}^c} \cup \bo')$ with $(k_1,m_1)=(k,m)$ and $|k_s|,|m_s| > n$ then there is a path in $\bo_{{\L_n}^c} \cup \bo'$ connecting $\bm{x}_{k,m} \in \bo'_{C(k,m)}$ to the point in $\bo_{C(k_s,m_s)}$.
    \qed
\end{proofsect}

\subsection{Percolation of $\mathcal{L}_{X_n}$.}\label{sec-percolation}
It remains to show that the assumptions of Corollary \ref{mixed-site bond positive percolation probability} are satisfied for the measures $(\mathcal{L}_{X_n})_{n=1}^\infty.$ The measures $\mathcal{L}_{X_n}$ can be considered as measures on $\{0,1\}^{\Z^2 \cup \mathcal{B}}$ where all bonds are opened with probability 1. For any $(k,m) \in \{-n,...,n\}^2$ and $X \in \{0,1\}^{\Z^2}$ satisfying $X(i,j) = 1$ for all $(i,j) \notin \{-n,...,n\}^2,$ 
\begin{align}
    &    \Csf_{\L_n,\omega}^{\ssup{\text{site}}}    \bigg(X_n(k,m) = 1 \bigg| X_n(i,j) = X(i,j) \text{ for } (i,j) \neq (k,m) \bigg) \label{k,m open conditioned on everything else} \\
    = \ &  \Csf_{\L_n,\omega}^{\ssup{\text{site}}}   \bigg(  \Csf_{\L_n,\omega}^{\ssup{\text{site}}} (F_{k,m} \cap O_{k,m}| \bF_{{C(k,m)}^c}) \bigg| X_n(i,j) = X(i,j) \text{ for } (i,j) \neq (k,m) \bigg) \nonumber\\
    \geq \ &\essinf_{\bo' \in \bO_{C(k,m)^c}}  \Csf_{\L_n,\omega}^{\ssup{\text{site}}}  \bigg(F_{k,m} \cap O_{k,m} \bigg| \bm{\pr}_{C(k,m)^c}  = \bo' \bigg),
    \label{essential infimum which needs lower bound}
\end{align}
where the essential infimum is taken with respect to $  \Csf_{\L_n,\omega}^{\ssup{\text{site}}}  \circ {\bm{\pr}_{C(k,m)^{\rm c}}}^{-1}.$ Since we are dealing with standard Borel spaces, the regular conditional probability in (\ref{essential infimum which needs lower bound}) is guaranteed to exist. For $(k,m) \notin \{-n,...,n\}^2$, there is only one point in each cell, so the expression (\ref{k,m open conditioned on everything else}) is equal to $\hat{p}$, which is in turn greater than (\ref{essential infimum which needs lower bound}) since the latter is at most $\hat{p}^{64}$. Therefore it is sufficient to show that (\ref{essential infimum which needs lower bound}) is greater than the critical probability for site percolation on $\Z^2$, denoted $p^{\text{site}}_{\rm c}(\Z^2)$. 

First we will bound the probability of the event $F_{k,m}$ from below. For $\Delta \subset \L_n$ the regular conditional probability $$\O_{\Delta^{\rm c}} \times \Fcal \ni (\omega',B) \mapsto M_{\L_n,\omega} (B | \pr_{\Delta^{\rm c}}  = \omega')$$ is given $M_{\L_n,\omega} \circ \pr_{\Delta^{\rm c}}^{-1}$-almost everywhere by the function  
\begin{align*}
    (\omega',B) \mapsto \frac{\int \1_{B}(\omega' \cup \omega'') h_{\L_n,\omega}(\omega' \cup \omega'') \,\ex^{H_{\L_n,\omega}(\omega'' \cup \omega'_{\L_n})} \,\Pi^z_{\Delta}(\d \omega'')}{\int h_{\L_n,\omega}(\omega \cup \omega') \,\ex^{H_{\L_n,\omega}(\omega')} \Pi^z_{\L_n}(\d \omega')}.
\end{align*}
The proof of this follows that of the analogous case in \cite[page 40-41]{E14}, or \cite{AE16}. For the rest of this section let $\epsilon = \frac{1}{2}(1 - p_{\rm c}^{\text{site}}(\Z^2)).$

\medskip

\begin{figure}
	\centering
	\includegraphics[width=0.7\linewidth]{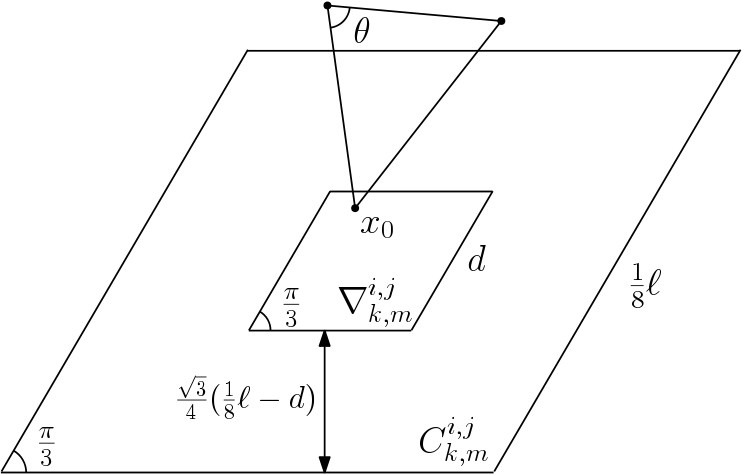}
	\caption{}
	\label{figure - inserting point into a sub cell}	
\end{figure}

\begin{lemma}
\label{hardcore at least one particle}
	Suppose $\alpha_0 < \sin^{-1} \left( \frac{3}{64} \right)$, $64r < 3R,$ \\ $\ell \in \left( \frac{64}{\sqrt{3}}(r \ \vee \ R\sin(\alpha_0)), \sqrt{3}R \right)$ and 
	\begin{align}
	    \label{z0 definition}
		z > z_0''(\ell,q,r,R,\alpha_0) := \frac{64q^{\frac{2\pi}{\alpha_0} - 1}}{\epsilon \left( \frac{\ell}{8} - \frac{8}{\sqrt{3}}(r \vee R\sin(\alpha_0)) \right)^2}.
	\end{align}
	Then for any pseudo-periodic boundary condition $\omega \in \Gamma$ and any sub-cell $C^{i,j}_{k,m}$ with $|k|,|m| \leq n$,
	\begin{align*}
		M_{\L_n,\omega}(N_{{C^{i,j}_{k,m}}} \geq 1 | \pr_{({C^{i,j}_{k,m}})^{\rm c}} = \omega') > 1 - \frac{\epsilon}{64}.
	\end{align*}
	for $M_{\L_n,\omega} \circ \pr_{(C^{i,j}_{k,m})^{\rm c}}^{-1}$-almost all $\omega'.$
\end{lemma}
\begin{proofsect}{Proof}
    Assume $M_{\L_n,\omega}(N_{{C^{i,j}_{k,m}}} = 0 | \pr_{({C^{i,j}_{k,m}})^{\rm c}} = \omega') > 0,$ else the result is trivial. This implies that
    \begin{align}
        \label{hardcore condition is satisfied if no points are added}
        \Psi(\tau) = 0 \text{ for all } \tau \in \Del_{3,\L_n}(\omega').
    \end{align}
    Define $\nabla_{k,m}^{i,j}$ to be the rhombus of side length $d = \frac{\ell}{8} - \frac{8}{\sqrt{3}}(r \vee R\sin(\alpha_0))$ which is a contraction of $C^{i,j}_{k,m}$ about its center point (see figure \ref{figure - inserting point into a sub cell}). We first claim that for $M_{\L_n,\omega} \circ \pr_{(C^{i,j}_{k,m})^{\rm c}}^{-1}$-almost all $\omega'$ and $x \in \nabla_{k,m}^{i,j}$,
    \begin{align}
        \label{hardcore condition satisfied if point is added in nabla}
        \Psi(\tau) = 0 \text{ for all } \tau \in \Del_{3,\L_n}(\omega' \cup x).  
    \end{align}
    It suffices to only consider the triangles $\tau \in T^+_{x,\omega'}.$ Any edge $\{x,y\} \in \Del_2(\omega' \cup x)$ must satisfy 
    $$|x-y|> \frac{\sqrt{3}}{4}(\frac{1}{8}\ell - d) = 2(r \vee R \sin(\alpha_0)) > 2r\,.
    $$ Thus $\delta(\tau) > r$ for all $\tau \in T^+_{x,\omega'}.$ By Lemma \ref{point insertion} and (\ref{hardcore condition is satisfied if no points are added}) we also have $\delta(\tau) <R.$

    By the same argument used to compare the angles $\theta_1$ and $\theta_2$ in Lemma \ref{point insertion}, if $\theta$ is an angle belonging to a triangle $\tau \in T^+_{x,\omega'}$ which is subtended at $x$ then $\theta \geq \alpha_0$. All edges $\{x,y\} \in \Del_2(\omega' \cup x)$ must have length at least $\frac{\sqrt{3}}{4}(\frac{1}{8}\ell - d) \geq 2R\sin(\alpha_0)$, so by the law of sines if $\tau \in T^+_{x,\omega'}$ and $\theta$ is an angle of $\tau$ not subtended at $x_0$ then 
	\begin{align*}
	    \sin(\theta) \geq \frac{2R \sin(\alpha_0)}{2\delta(\tau)} \geq \sin(\alpha_0),
	\end{align*} 
	and so $\theta \geq \alpha_0$ since we know that $\alpha_0 \leq \pi/3$. This completes the proof of (\ref{hardcore condition satisfied if point is added in nabla}). Together with (\ref{hardcore condition is satisfied if no points are added}) this implies that
	\begin{align}
	    \label{Hamiltonians are 0 whether adding point or not} 
	    \ex^{-H_{\L_n,\omega} \left(\omega_{\L_n}' \cup x \right)} = \ex^{-H_{\L_n,\omega}\left(\omega'_{\L_n}\right)} = 1.
	\end{align}
	For the second half of the proof we compute the lower bound by applying the formula 
	\begin{align*}
		\int f(\omega'') \, \Pi^z_{\nabla_{k,m}^{i,j}}(\d \omega'') = \ex^{-z|\nabla_{k,m}^{i,j}|}\sum_{n=0}^\infty \frac{z^n}{n!} \int_{(\nabla_{k,m}^{i,j})^n} f(\{ x_1, ... , x_n \}) \,\d x_1, .. , \d x_n
	\end{align*} 
	(valid for bounded measurable functions $f:\O_{\nabla_{k,m}^{i,j}} \rightarrow [0,\infty)$), restricting the domain of integration and applying (\ref{Hamiltonians are 0 whether adding point or not}) and Lemma \ref{hardcore_papangelou_intensity}. The computation is as follows:
	\begin{align*}
		& \hspace{0.5cm} \frac{M_{\L_n,\omega}(N_{{C^{i,j}_{k,m}}} = 1 | \pr_{({C^{i,j}_{k,m}})^{\rm c}} = \omega')}{M_{\L_n,\omega}(N_{{C^{i,j}_{k,m}}} = 0 | \pr_{({C^{i,j}_{k,m}})^{\rm c}} = \omega')} \\[0.5em]
		&= \frac{\int \1_{\{N_{{C^{i,j}_{k,m}}} = 1\}}(\omega'') h_{\L_n,\omega}(\omega' \cup \omega'') \ex^{-H_{\L_n,\omega}(\omega'_{\L_n} \cup \omega'')} \,\Pi^z_{C^{i,j}_{k,m}}(\d \omega'')}{\int \1_{\{N_{{C^{i,j}_{k,m}}} = 0\}}(\omega'') h_{\L_n,\omega}(\omega')\, \ex^{-H_{\L_n,\omega}(\omega'_{\L_n})} \,\Pi^z_{C^{i,j}_{k,m}}(\d \omega'')} \\[0.5em]
		&= \frac{z \ex^{-z|C^{i,j}_{k,m}|} \int_{C^{i,j}_{k,m}} h_{\L_n,\omega}(\omega' \cup x) \,\ex^{-H_{\L_n,\omega}(\omega'_{\L_n} \cup x)} \,\d x }{\ex ^{-z|C^{i,j}_{k,m}|}  h_{\L_n,\omega}(\omega')}  \\[0.5em]
		&\geq z  \int_{\nabla^{i,j}_{k,m}} \frac{h_{\L_n,\omega}(\omega' \cup x)}{h_{\L_n,\omega}(\omega')} \,\d x \\[0.5em]
		&\geq z q^{1-\frac{2\pi}{\alpha_0}} |\nabla^{i,j}_{k,m}| \\[0.5em]
		&= z q^{1-\frac{2\pi}{\alpha_0}} \left( \frac{\ell}{8} - \frac{8}{\sqrt{3}}(r \vee R\sin(\alpha_0)) \right)^2.
	\end{align*}
	Finally,
	\begin{align*}
		M_{\L_n,\omega}(N_{{C^{i,j}_{k,m}}} \geq 1 | \pr_{({C^{i,j}_{k,m}})^{\rm c}} = \omega') 
		&\geq 1 - \frac{M_{\L_n,\omega}(N_{{C^{i,j}_{k,m}}} = 1 | \pr_{({C^{i,j}_{k,m}})^{\rm c}} = \omega')}{M_{\L_n,\omega}(N_{{C^{i,j}_{k,m}}} = 0 | \pr_{({C^{i,j}_{k,m}})^{\rm c}} = \omega')} \\[0.5em]
		&\geq 1 - \frac{q^{\frac{2\pi}{\alpha_0} - 1}}{z \left( \frac{\ell}{8} - \frac{8}{\sqrt{3}}(r \vee R\sin(\alpha_0)) \right)^2} \\[0.5em]
		&> 1 - \frac{\epsilon}{64}
	\end{align*}
	by assumption (\ref{z0 definition}).
	\qed
\end{proofsect}

\begin{cor}
    \label{lower bound on event Fkm}
    If the assumptions of Lemma \ref{hardcore at least one particle} are satisfied, then 
    \begin{align*}
		M_{\L_n,\omega}(F_{k,m}| \pr_{C(k,m)^{\rm c}} = \omega') > 1 - \epsilon.
	\end{align*}
	for $M_{\L_n,\omega} \circ \pr_{C(k,m)^{\rm c}}$-almost all $\omega'.$
\end{cor}
\begin{proofsect}{Proof}
    \begin{align*}
        &M_{\L_n,\omega}(F_{k,m}| \pr_{C(k,m)^{\rm c}} = \omega') \\
        \geq \ &1 - \sum_{0 \leq i,j \leq 7}  M_{\L_n,\omega}(N_{{C^{i,j}_{k,m}}} = 0 | \pr_{C(k,m)^{\rm c}} = \omega') \\
        \geq \ &1 - \sum_{0 \leq i,j \leq 7} \int M_{\L_n,\omega}(N_{{C^{i,j}_{k,m}}} = 0 | \pr_{(C_{k,m}^{i,j})^{\rm c}} = \omega') M_{\L_n,\omega} \circ \pr_{(C^{i,j}_{k,m})^{\rm c}}^{-1}(\d \omega') \\
        > \ &1 - \epsilon. \qedhere
    \end{align*}
    \qed
\end{proofsect}

The final component we need to finish the proof is an upper bound on the number of particles in a cell $C(k,l)$. If $H_{\L_n,\omega}(\omega') < \infty$ then $\{x,y\} \in \Del_{2,\L_n}(\omega_{\L_n^c} \cup \omega') \implies |x-y| \geq 2r\sin(\alpha_0)$ by the law of sines. Since the Delaunay graph is a nearest neighbour graph, this means that no two particles are within a distance of $2r\sin(\alpha_0)$ of one another. Therefore 
\begin{align}
    \label{upper bound m definition}
	m(\ell,r,\alpha_0) := \left( \frac{\ell + 2r\sin(\alpha_0)}{r\sin(\alpha_0)} \right)^2 = \left( \frac{\ell}{r\sin(\alpha_0)} + 2 \right)^2,
\end{align}
which is an upper bound on the number of non-overlapping circles with radius $r\sin\alpha_0$ that can fit inside a rhombus with side length $\ell + 2r\sin\alpha_0$, is an upper bound for $|\omega'|.$ 

We can now prove the existence of the required lower bound on (\ref{essential infimum which needs lower bound}).
\begin{proposition}
    \label{lower bound on Csite conditional probabilities}
	Suppose $\alpha_0 < \sin^{-1}(\frac{3}{64})$, $64r < 3R$ and $\ell \in ( \frac{64}{\sqrt{3}}(r \ \vee \ R\sin(\alpha_0)), \sqrt{3}R)$. If  
	$$
		\beta > \beta'_0(\ell,q,R,r,\alpha_0) 
    := \ \frac{\frac{3\sqrt{3}}{4}q^2 R^2 }{(1-\epsilon)^{-1/m(\ell,r,\alpha_0)} - 1}
	$$
	and
	$
		z > z_0''(\ell,q,R,r,\alpha_0)
	$
	then there exists $c>0$ such that for any $n \in \N$, any $|k|,|m| \leq n$ and any pseudo-periodic boundary condition $\omega \in \Gamma$, 
	\begin{align*}
		   \Csf_{\L_n,\omega}^{\ssup{\rm site}} \bigg(F_{k,m} \cap O_{k,m} \bigg| \bm{\pr}_{C(k,m)^{\rm c}}  = \bo' \bigg) \geq c > 0
	\end{align*}
	for $\Csf^{\ssup{\text{site}}}_{\L_n,\omega} \circ \bm{\pr}_{C(k,m)^{\rm c}}^{-1}$-almost all $\bo'$.
\end{proposition}
\begin{proofsect}{Proof}
    Using the lower bound from Corollary \ref{lower bound on event Fkm} and the lower bound $\beta_0$ we have
    \begin{align*}
	\Csf_{\L_n,\omega}^{\ssup{\text{site}}}	    (F_{k,m} \cap O_{k,m} | \bm{\pr}_{C(k,m)^{\rm c}}  = \bo' )
		&\geq \hat{p}^m M_{\L_n,\omega} ( F_{k,m} | \pr_{C(k,m)^{\rm c}} = \omega' ) \\
		&\geq \left( \frac{1}{\frac{3\sqrt{3}}{4\beta}q^2 R^2 + 1} \right)^m (1-\epsilon) \\
		&\geq \left(\frac{1}{(1-\epsilon)^{-\frac{1}{m}}} \right)^m (1-\epsilon) \\
		&= (1-\epsilon)^2 \\
		&> 1 - 2\epsilon = p_{\rm c}^{\rm site}(\Z^2) > 0. \qedhere
	\end{align*}
	\qed
\end{proofsect}

Recall that at the beginning of section \ref{coarse graining section} we assumed that the parameters $\rho,\ell$ and $zq$ satisfied the requirements of Proposition \ref{existence theorem hardcore}. We need to check that these requirements can be satisfied simultaneously with those of the previous proposition. The only possible conflict relates to the parameter $\ell$. We require that $\ell \in \left(\frac{64}{\sqrt{3}}(r \ \vee \ R\sin(\alpha_0)), \sqrt{3}R \right) \cap \left( \frac{r}{L(\rho)}, \frac{R}{U(\rho)} \right):= I_0(R,r,\alpha_0,\rho).$ This set is non-empty for small enough $\rho$ since 
\begin{align*}
    \frac{64}{\sqrt{3}}(r \ \vee \ R\sin(\alpha_0)) < \sqrt{3}R = \lim_{\rho \rightarrow 0} \frac{R}{U(\rho)}
\end{align*}
and 
\begin{align*}
    \sqrt{3} R > \sqrt{3}r = \lim_{\rho \rightarrow 0} \frac{r}{L(\rho)}.
\end{align*}
We therefore have the following, from which Proposition~\ref{Prop-per} and Theorem~\ref{THM-main} is derived by selecting particular values of $\rho$ and $\ell.$
\begin{cor}
    There exists $\rho'_0(R,r,\alpha_0)>0$ such that $I_0 \neq \emptyset$ if $\rho < \rho'_0$. Moreover, if
    \begin{align*}
        \alpha_0 &< \sin^{-1}(3/64), \\
        64r &< 3R, \\
        \rho &< \rho'_0(R,r,\alpha_0), \\
        \ell &\in I_0(R,r,\alpha_0,\rho), \\
        \beta &> \beta_0(\ell,q,R,r,\alpha_0), \\
        \text{ and } \qquad z &> z_0''(\ell,q,R,r,\alpha_0) \ \vee \ (1/q) z'_0(\beta, \rho, \ell)
    \end{align*}
    then there exists at least $q$ translation-invariant Delaunay continuum Potts measures for $\bDel_3, z$ and $\varphi.$
\end{cor}
{\color{blue} 
}

\section*{Appendix}
\begin{appendices}

 \section{Pseudo-periodic configurations}\label{pseudoperiodic}
  We define pseudo-periodic configurations as in \cite{DDG12}. We first obtain a partition of $ \R^2 $  into rhombuses. Pick a length scale $ \ell>0 $ and consider the matrix $$ M=\left(\begin{matrix}M_1 & M_2 \end{matrix}\right) =\left(\begin{matrix} \ell & \ell/2\\0 &\sqrt{3}/2 \ell\end{matrix}\right).$$
 Note that $ |M_i|=\ell, i=1,2 $, and $ \angle(M_1,M_2)=\pi/3 $. For each $ (k,l)\in\Z^2 $ we define the cell
 \begin{equation}\label{cell}
 \Delta_{k,l}=\{Mx\in\R^2\colon x-(k,l)\in[-1/2,1/2)^2\}
 \end{equation} with area $ |\Delta_{k,l}|=\frac{\sqrt{3}}{2}\ell^2 $. For example, $ \Delta_{0,0} $ is the rhombus with corners $ (3\ell/4,\sqrt{3}\ell/4), (\ell/4,-\sqrt{3}\ell/4),(-3\ell/4,-\sqrt{3}\ell/4),(-\ell/4,\sqrt{3}\ell/4) $,  and horizontal side length of $ \ell $. These cells constitute a periodic partition of $ \R^2 $ into rhombuses. Let 
 \begin{equation}\label{ppc}
\Gamma=\{\om\in\O\colon\theta_{Mz}(\om_{\Delta_{k,l}})\in B\mbox{ for all } z=(k,l)\in\Z^2, B \mbox{ measurable set of } \O_{\Delta_{0,0}}\setminus\{\emptyset\}  \} 
 \end{equation} be the set of all configurations whose restriction to an arbitrary cell, when shifted back to $ \Delta_{0,0} $, belongs to the  measurable set $ B $ for all measurable sets $ B$ of $\O_{\Delta_{0,0}}\setminus\{\emptyset\} $. Elements of $ \Gamma $ are called \textbf{pseudo-periodic} configurations. We define marked pseudo-periodic configurations in an analogous way.

\section{Topology of local convergence}\label{AppC}
We write $ \Mcal_1^{\theta}(\bO) $ (resp. $ \Mcal_1^{\theta}(\O) $) for the set of all shift-invariant probability measures on $ (\bO,\boldsymbol{\Fcal}) $ (resp. $ (\O,\Fcal) $). A measurable function $ f\colon\bO\to\R $ is called local and tame if
$$
f(\bo)=f(\bo_\L)\quad\mbox{ and }\quad |f(\bo)|\le a N_\L(\bo)+b
$$ for all $ \bo\in\bO $ and some $ \L\Subset\R^2 $ and suitable constants $ a,b\ge 0 $. Let $ \Lscr $ be the set of all local and tame functions. The topology of local convergence, or $\Lscr$-topology, on $ \Mcal_1^{\theta}(\bO) $ is then defined as the weak$*$ topology induced by $ \Lscr $, i.e., as the smallest topology for which the mappings $ P\mapsto \int f\d P $ with $ f\in\Lscr $ are continuous.

\section{Mixed site-bond percolation}\label{AppD}

Here we take a small detour to prove a technical result about mixed site-bond percolation on $\Z^2$ which will prove useful when carrying out the aforementioned coarse-graining procedure. This result is not strictly necessary, and in fact \cite[Lemma 1]{Rus82} would suffice for our purposes since we only need to consider site percolation. We choose to include this result here because there are situations (for instance the model considered in \cite{AE19}) when one needs to use mixed site-bond percolation to accomplish the coarse-graining argument. The result states that if all conditional probabilities are uniformly bounded from below then the percolation probability is greater than the percolation probability of Bernoulli site-bond percolation. For more results regarding comparisons of site percolation measures with product measures see \cite{LSS97}.

Let $\mathcal{B} = \{ \{x,y\} \subset \Z^2 : |x-y| = 1 \}$ denote the set of edges (or \textit{bonds}) between neighbouring vertices in $\Z^2$, and $\O := \{0,1\}^{\Z^2 \cup \mathcal{B}}$. For $\om \in \O$ we say that a site or bond $x \in \Z^2 \cup \mathcal{B}$ is \textit{open} if $\om(x)=1$ and \textit{closed} otherwise. We will use the shorthand $\om_x$ and $\om_{\{x,y\}}$ in place of $\om(x)$ and $\om(\{x,y\})$ respectively. Let
\begin{align*}
	C(\omega) :=
	\left\{ x \in \Z^2 \ 
    \begin{tabular}{|c}
    	$\text{There exists a path } 0=x_1,...,x_n=x \text{ s.t } \om_{x_i}=1 $\\
    	$\text{ for } i \in [n] \text { and } \om_{ \{ x_i,x_{i+1} \} } =1 \text{ for } i\in [n-1]. $
	\end{tabular}
	\right\}
\end{align*}
This is the \textit{open cluster} around the origin. Similarly, if $\om \in \O$ or $\om \in \tilde{\O}:= \{0,1\}^{\Z^2}$ we define the \textit{open site cluster} as follows:
\begin{align*}
	C_s(\omega) :=
	\left\{ x \in \Z^2 \ 
    \begin{tabular}{|c}
    	$\text{There exists a path } 0=x_0,,..,x_n=x$ \\
    	$\text{ such that } \om_{x_i}=1 \text{ for } i \in [n] $ \\
	\end{tabular}
	\right\}.
\end{align*}
The event that $0$ is connected to the set $A$ is
\begin{align*}
	\{ 0 \leftrightarrow A \} = \{ \om \in \O \ | \ \text{There exists } x \in A \cap C(\om) \},
\end{align*}
and the event that percolation occurs is
\begin{align*}
	\{ 0 \leftrightarrow \infty\} := \{ \om \in \O \ | \ |C(\om)| = \infty\} = \bigcap_{n \in \N} \{ 0 \leftrightarrow \L_n^c\},
\end{align*}
where $\L_n = [-n,n]^2 \cap \Z^2$.

Let $\mu_{p,p'}$ denote the measure for which each site is opened independently with probability $p$ and each bond is opened independently with probability $p'$.


\begin{lemma}\label{percolation greater than bernoulli}
	If $\P$ is a measure on $\O$ satisfying the conditions: 
	\begin{enumerate}
		\item 
			For all $x \in \Z^2$ and $\om' \in \O,$
			\begin{align} \label{open site greater than p}
	            \P
	            \left( \om_x = 1 \ 
                \begin{tabular}{|c}
    	            $\om_{\{w,z\}}  = \om'_{\{w,z\}} \text{ for all } \{w,z\}$ \\
    	            and $\om_y = \om'_y \text{ for all } y \neq x$ \\
	            \end{tabular}
	            \right) \geq p.
            \end{align}
		\item 
			For all $\{x,y\} \in \mathcal{B}$ and $\om' \in \O$ satisfying $\om'_x = \om'_y = 1,$
			\begin{align}\label{open bond greater than p'}
	            \P
	            \left( \om_{\{x,y\}} = 1 \ 
                \begin{tabular}{|c}
    	            $\om_{\{w,z\}}  = \om'_{\{w,z\}} \text{ for } \{w,z\} \neq \{x,y\}$ \\
    	            and $\om_y = \om'_y \text{ for all } y $ \\
	            \end{tabular}
	            \right) \geq p'.
            \end{align}
	\end{enumerate}
	Then $\P(0 \leftrightarrow \infty) \geq \mu_{p,p'}(0 \leftrightarrow \infty).$
\end{lemma}
\begin{proofsect}{Proof}
	Let $\tilde{\P}$ and $\tilde{\mu}_{p,p'}$ denote the marginal measures of $\P$ and $\mu_{p,p'}$ on $\tilde{\O}$. Inequality (\ref{open site greater than p}) will allow us to couple these measures together. We start by identifying $\Z^2$ with the natural numbers via an arbitrary ordering. Let $E_k^+ \ (E_k^-)$ be the event that the site $k$ is open (closed) and let $\om^{(k)} = \{ \om' \in \tilde{\O} \ | \ \om'_i = \om_i \text{ for all } i < k \}$. We define the measure $m$ on $\tilde{\O} \times \tilde{\O}$ inductively by setting (as was done in \cite{Rus82})
	\begin{align*}
		& m(E_1^+ \times E_1^+) = p, && m(E_1^+ \times E_1^-) = \P(E_1^+) - p, \\
		& m(E_1^- \times E_1^+) = 0, && m(E_1^- \times E_1^-) = 1-\P(E_1^+) ,
	\end{align*}
	and then for $k \geq 2$ and any $\zeta, \omega \in \tilde{\O},$
	\begin{align*}
		& m(E_k^+ \times E_k^+|\zeta^{(k)} \times \omega^{(k)}) = p, && m(E_k^+ \times E_k^-|\zeta^{(k)} \times \omega^{(k)}) = \P(E_k^+ \ | \ \omega^{(k)}) - p, \\
		& m(E_k^- \times E_k^+|\zeta^{(k)} \times \omega^{(k)}) = 0, && m(E_k^- \times E_k^-|\zeta^{(k)} \times \omega^{(k)}) = 1 - \P(E_k^+ \ | \ \omega^{(k)}).
	\end{align*}
	This measure satisfies the following:
	\begin{align*}
		\P(0 \leftrightarrow \Lambda_n^{\rm c}) 
		&= \int \P(0 \leftrightarrow \L_n^{\rm c} \ | \ C_s = C_s(\om))\ \tilde{\P}(\d \om) \nonumber\\
		&= \int \P(0 \leftrightarrow \Lambda_n^{\rm c} \ | \ C_s = C_s(\om))\ m(\d \om,\d \om'),
	\end{align*}
	and similarly 
	\begin{align*} 
		\mu_{p,p'}(0 \leftrightarrow \L_n^{\rm c}) 
		&= \int \mu_{p,p'}(0 \leftrightarrow \L_n^{\rm c} \ | \ C_s = C_s(\om'))\ m(\d \om,\d \om').
	\end{align*}
	Therefore if
	\[
		\int \P(0 \leftrightarrow \L_n^{\rm c} \ | \ C_s = C_s(\om)) - \mu_{p,p'}(0 \leftrightarrow \L_n^{\rm c} \ | \ C_s = C_s(\om ')) \ m(\d \om ,\d \om ') \geq 0 
	\]
	then we have
	\[
		\P(0 \leftrightarrow \L_n^{\rm c}) \geq \mu_{p,p'}(0 \leftrightarrow \L_n^c),
	\]
	and so to complete the proof we need only to show the former inequality. First we note that if $\om \geq \om'$, then we have 
	\begin{align} \label{bernoulli Cs monotonicity}
		\mu_{p,p'}(0 \leftrightarrow \L_n^{\rm c} \ | \ C_s = C_s(\om)) \geq \mu_{p,p'}(0 \leftrightarrow \L_n^c \ | \ C_s = C_s(\om')), 
	\end{align} 
	since the event $0 \leftrightarrow \L_n^{\rm c}$ only depends on the bonds between sites in $C_s$. Secondly, since all sites in $C_s$ are open, we can conclude using property (\ref{open bond greater than p'}) that for all $\om \in \tilde{\O}$
	\begin{align} \label{bond coupling}
		\P(0 \leftrightarrow \L_n^{\rm c} \ | \ C_s = C_s(\om)) \geq \mu_{p,p'}(0 \leftrightarrow \L_n^{\rm c} \ | \ C_s = C_s(\om)).
	\end{align}
	This can be shown by another coupling on $\{0,1\}^{B(\om)}$, where $B(\om)$ is the set of bonds between sites in $C_s(\om)$. Since $\om \geq \om'$ almost surely with respect to $m$ we can now conclude using (\ref{bernoulli Cs monotonicity}) and (\ref{bond coupling}) that: 
	\begin{align*}
		&\int \P(0 \leftrightarrow \L_n^{\rm c} \ | \ C_s = C_s(\om))\ - \mu_{p,p'}(0 \leftrightarrow \L_n^{\rm c} \ | \ C_s = C_s(\om'm)) \ m(\d \om,\d \om') \\
		\geq &\int \P(0 \leftrightarrow \L_n^{\rm c} \ | \ C_s = C_s(\om))\ - \mu_{p,p'}(0 \leftrightarrow \L_n^{\rm c} \ | \ C_s = C_s(\om)) \ m(\d \om,\d \om') \\
		\geq &\ 0. \qedhere
	\end{align*}
\end{proofsect}

\begin{cor}
    \label{mixed-site bond positive percolation probability}
	Let $p_{\rm c}$ denote the critical probability for site percolation on $\Z^2$. If $\P$ is a measure on $\O$ such that
	\begin{enumerate}
		\item 
			For all $x \in \Z^2$ and $\om' \in \O,$
			\begin{align*} 
	            \P
	            \left( \om_x = 1 \ 
                \begin{tabular}{|c}
    	            $\om_{\{w,z\}}  = \om'_{\{w,z\}} \text{ for all } \{w,z\}$ \\
    	            and $\om_y = \om'_y \text{ for all } y \neq x$ \\
	            \end{tabular}
	            \right) \geq \sqrt{p_{\rm c}}.
            \end{align*}
		\item 
			For all $\{ x,y \} \in \mathcal{B}$ and all $\om' \in \O$ satisfying $\om'_x = \om'_y = 1,$
			\begin{align*}
	            \P
	            \left( \om_{\{x,y\}} = 1 \ 
                \begin{tabular}{|c}
    	            $\om_{\{w,z\}}  = \om'_{\{w,z\}} \text{ for } \{w,z\} \neq \{x,y\}$ \\
    	            and $\om_y = \om'_y \text{ for all } y $ \\
	            \end{tabular}
	            \right) \geq \sqrt{p_{\rm c}}.
            \end{align*}
	\end{enumerate}
	Then $\P(0 \leftrightarrow \infty) > 0.$
\end{cor}

\begin{proofsect}{Proof}
	By Lemma (\ref{percolation greater than bernoulli}) there exists $\epsilon>0$ such that
	\begin{align*}
		\P(0 \leftrightarrow \infty) \geq \mu_{\sqrt{p_{\rm c}}+\epsilon,\sqrt{p_{\rm c}}+\epsilon}(0 \leftrightarrow \infty).
	\end{align*}
	By applying inequality (4) of \cite{Ha80} we can see that the right hand side is greater than $\mu_{(\sqrt{p_{\rm c}}+\epsilon)^2,1}(0 \leftrightarrow \infty)$, which is greater than $0$ since $(\sqrt{p_{\rm c}}+\epsilon)^2 > p_{\rm c}$.

\qed
\end{proofsect}

\end{appendices}
\section*{Acknowledgments}
Shannon Horrigan thanks MASDOC - the centre of doctoral training (CDT) at Warwick.


\bibliographystyle{Martin}

\bibliography{Adams_Continuum}

\end{document}